\documentclass[11pt]{article}
\usepackage[a4paper,margin=26mm]{geometry}
\usepackage[T1]{fontenc}
\usepackage{lmodern}
\usepackage{amsmath,amssymb,amsthm,mathtools}
\usepackage{booktabs,tabularx,array}
\usepackage[hidelinks,unicode]{hyperref}
\usepackage{microtype}
\newtheorem{theorem}{Theorem}[section]
\newtheorem{lemma}[theorem]{Lemma}
\newtheorem{proposition}[theorem]{Proposition}
\newtheorem{corollary}[theorem]{Corollary}
\theoremstyle{definition}

\theoremstyle{remark}
\newtheorem{remark}[theorem]{Remark}
\numberwithin{equation}{section}
\newcommand{\Id}{\mathcal I}
\newcommand{\Proj}{\mathcal P}
\newcommand{\B}{\mathcal B}
\newcommand{\tr}{\operatorname{Tr}}
\newcommand{\norm}[1]{\left\|#1\right\|}

\newcommand{\C}{\mathbb C}
\newcommand{\R}{\mathbb R}
\newcommand{\N}{\mathbb N}
\DeclareMathOperator{\Ran}{Ran}
\DeclareMathOperator{\rank}{rank}

\hypersetup{pdftitle={Lp-isometries on Idempotents},pdfsubject={Fixed-rank rigidity and Lp perturbation classes}}
\title{$L^p$-isometries on idempotents\\[4pt]
\large Fixed-rank rigidity, the operator-norm case,\\and perturbation classes}
\author{Wenhua Qian$^{1, \dag}$, Yuze Ren$^1$, Wenming Wu$^1$}
\date{}
\begin{document}
\maketitle
\begingroup
\renewcommand{\thefootnote}{}
\footnotetext{Qian ($^\dag$: Corresponding Author) is supported in part by the Natural Science Foundation of Chongqing Science and Technology Commission (CSTB2025NSCQ-GPX1002), and Wu is supported in part by NSF of China (No.\ 12571130) and the Natural Science Foundation of Chongqing Science and Technology Commission (Grant No.\ CSTB2025NSCQ-LZX0058).}
\endgroup
\begin{abstract}
We study surjective $L^p$-isometries on spaces of bounded idempotents in an infinite-dimensional complex Hilbert space. On each fixed positive finite-rank stratum, every surjective $L^p$-isometry, $1\leq p<\infty$, is implemented by a unitary or antiunitary operator, possibly composed with adjunction; no separability is required. The same conclusion holds for the operator norm, viewed as the $L^\infty$-norm. On a separable space, operator-norm isometries of the infinite-rank, infinite-corank sector admit the additional complement operation, as do surjective $L^p$-isometries between individual $L^p$-perturbation classes for $1\leq p<\infty$. For the extended $L^p$-metric on all idempotents, implementations may be chosen independently on finite-distance classes, subject to an explicit compatibility condition for matching infinite-rank classes. The proofs use similarity coordinates and metric blow-ups, followed by tangent-space splittings, nearest projections at the $L^1$ endpoint, or sharp Fredholm-index distance estimates. We also establish uniqueness of standard implementations and use essential-spectral invariants to obtain continuum many perturbation-class isometry types.
\end{abstract}

\tableofcontents

\section{\texorpdfstring{Introduction and the $L^p$-metric setting}{Introduction and the Lp-metric setting}}\label{sec:introduction}
Let $H$ be an infinite-dimensional complex Hilbert space. An idempotent $E\in\B(H)$ satisfies $E^2=E$; it is an orthogonal projection precisely when $E=E^*$. Every bounded idempotent determines a decomposition $H=\Ran E\dotplus\ker E$ into complementary closed subspaces, with $E$ acting as the projection onto the first summand along the second. Thus an orthogonal projection is determined by its range, whereas a general idempotent also records the choice of a complementary subspace. Even at a fixed positive finite rank, the set of idempotents is nonlinear and unbounded. These additional degrees of freedom make its metric geometry substantially richer than that of the corresponding Grassmann space.

A central problem in the geometry of operator spaces is to determine how much algebraic structure can be recovered from distances on a nonlinear subset. Wigner's theorem provides a fundamental example: a bijection of rank-one orthogonal projections preserving transition probabilities is implemented by a unitary or antiunitary operator \cite{ref17}. For rank-one projections, transition probabilities also determine distances in the operator norm and in every $L^p$-norm considered here. This connects quantum symmetries with the metric rigidity of Grassmann spaces. In contrast with isometries of whole real normed spaces, to which the Mazur--Ulam theorem applies \cite{ref11}, an isometry of an idempotent space has no immediate affine extension.

The classification of Grassmann-space isometries has developed in several stages. Botelho, Jamison and Moln\'ar \cite{ref1} classified surjective operator-norm isometries on fixed finite-rank projection spaces under a dimensional assumption. Geh\'er and \v Semrl \cite{ref5} removed that restriction, and subsequently treated the infinite-rank, infinite-corank projection space on a separable Hilbert space \cite{ref6}. For noncommutative $L^p$-metrics, Qian, Shen, Shi, Wu and Yuan \cite{ref13} classified surjective isometries between Grassmann spaces of fixed finite trace in semifinite factors, relating them to $*$-isomorphisms and $*$-anti-isomorphisms. These results establish a broad rigidity theory for orthogonal projections, but do not directly classify isometries on nonself-adjoint idempotents.

In the present type-I setting, for $1\leq p<\infty$ the noncommutative $L^p$-norms associated with the canonical trace on $\B(H)$ are the Schatten $p$-norms \cite{ref12,ref16}. They measure distances between idempotents through their differences, with $p=1$ giving the trace norm and $p=2$ the Hilbert--Schmidt norm; we also consider the operator norm, denoted by $p=\infty$. Differences within a fixed finite-rank stratum always have finite rank, so all these distances are finite there. For $1\leq p<\infty$, however, general infinite-rank idempotents can lie at infinite distance from one another, and the finite-distance relation singles out $L^p$-perturbation classes. The basic object throughout this paper is a surjective isometry on a specified stratum, sector, or perturbation class; the map is not assumed to be linear or to extend to an isometry of the ambient operator space.

For nonself-adjoint idempotents, Qian, Xiang, Wu and Yi \cite{ref14} classified surjective $L^p$-isometries on the rank-one stratum for $1<p<\infty$, obtaining unitary or antiunitary implementation, possibly composed with adjunction, through projection preservation. Qian, Shi, Wu, Xiang and Yuan \cite{ref15} established further operator-norm rigidity results, including the rank-one classification in Theorem~4.5 and the classification on the space of all nontrivial idempotents of a separable Hilbert space in Theorem~5.8. Their Remark~5.9 raised the fixed-rank problem for $k>1$. A classification on the full space does not by itself settle this problem: an isometry given only on one stratum is not assumed to extend to all idempotents. The same distinction applies to maps defined only on the infinite-rank, infinite-corank sector.

The questions addressed here are therefore whether these independently specified domains retain the rigidity of the projection spaces, and how that rigidity depends on $p$. We treat every positive finite rank for $1\leq p\leq\infty$ without separability, and the operator-norm infinite-rank, infinite-corank sector on separable $H$. For finite $p$ on separable $H$, we also study isometries between individual $L^p$-perturbation classes and the resulting extended metric on all idempotents. A perturbation class need not contain an orthogonal projection and may have several connected components. Consequently, its classification requires recovering structure from distances across the entire class, beyond the projection geometry available in a fixed-rank stratum.

The exposition begins with common local coordinates and a metric blow-up construction. The proofs then separate according to the geometry of the domain and the value of $p$. Except in the Hilbert--Schmidt case, tangent splittings detect orthogonal projections within a fixed-rank stratum. The index geometry of a finite-distance class instead detects orthogonal increments without first identifying its projection points. The $L^1$ endpoint and the $L^2$ case require their own arguments. The operator-norm cases are treated separately; they are not obtained by taking a formal limit as $p\to\infty$.

\subsection{\texorpdfstring{Idempotent spaces and $L^p$-distances}{Idempotent spaces and Lp-distances}}\label{intro:spaces}
Write
\[
\Id(H)=\{E\in\B(H):E^2=E\},\qquad
\Id_k(H)=\{E\in\Id(H):\rank E=k\},
\]
and let $\Proj_k(H)$ denote the self-adjoint members of $\Id_k(H)$. On a separable $H$, set
\[
\Id_{\infty,\infty}(H)
=\{E\in\Id(H):\dim\Ran E=\dim\ker E=\infty\},
\]
with projection subset $\Proj_{\infty,\infty}(H)$. The equality of the two infinite dimensions is used only in the separable setting; without separability, our results concern finite-rank strata.

Let $\tr$ be the canonical faithful normal semifinite trace on $\B(H)$. For $1\leq p<\infty$, we use the notation
\[
L^p(\B(H),\tr)=\{X\in\B(H):\tr(|X|^p)<\infty\},
\qquad \|X\|_p=\bigl(\tr(|X|^p)\bigr)^{1/p}.
\]
This is the type-I noncommutative $L^p$-space used throughout the paper; standard facts about these norms and their operator-ideal properties may be found in \cite{ref16}. We also write $L^\infty(\B(H),\tr)=\B(H)$ and $\|X\|_\infty=\|X\|$ for the operator norm. Unless explicitly marked by a subscript, a norm in the operator-norm sections means $\|\cdot\|_\infty$.

For $E,F\in\Id_k(H)$, the difference $E-F$ has rank at most $2k$. Hence
\[
d_p(E,F)=\|E-F\|_p
\]
is a finite metric on $\Id_k(H)$ for every $1\leq p\leq\infty$. A surjective $L^p$-isometry of this stratum is a bijection $\varphi$ satisfying
\[
d_p(\varphi(E),\varphi(F))=d_p(E,F)\qquad(E,F\in\Id_k(H)).
\]

For $1\leq p<\infty$, the full idempotent space carries the extended $L^p$-metric
\begin{equation}\label{intro:extendedmetric}
d_p(E,F)=
\begin{cases}
\|E-F\|_p,&E-F\in L^p(\B(H),\tr),\\
+\infty,&E-F\notin L^p(\B(H),\tr).
\end{cases}
\end{equation}
Its finite-distance classes are the $L^p$-perturbation classes
\[
\mathcal C_p(E_0)=\{E\in\Id(H):E-E_0\in L^p(\B(H),\tr)\}.
\]
The idempotents in such a class need not themselves belong to $L^p$; it is their differences that do. A finite-distance class may have several connected components, indexed below by a relative Fredholm index. We therefore distinguish a perturbation class from each of its connected components. No analogous finite-$p$ class decomposition is asserted for $p=\infty$.

\subsection{Standard implementations and the structure of the results}\label{intro:structure}
For a unitary or antiunitary operator $U$, the standard transformations are
\begin{equation}\label{intro:standard}
E\longmapsto UEU^*,\quad UE^*U^*,\quad U(I-E)U^*,\quad U(I-E^*)U^*.
\end{equation}
All four preserve the $L^p$-norm of differences for $1\leq p\leq\infty$. Only the first two preserve a positive finite-rank stratum on an infinite-dimensional space.

Section~\ref{sec:common} gives the common algebraic and metric tools. Sections~\ref{sec:operator-projections} and~\ref{sec:operator-classification} treat surjective $L^\infty$-isometries: finite rank requires no separability, while the infinite-rank, infinite-corank sector is considered on separable $H$. Section~\ref{sec:fixed-rank} proves the fixed-rank classification for $1\leq p<\infty$, again without separability. Section~\ref{sec:classes} treats individual $L^p$-perturbation classes and the extended $L^p$-metric on $\Id(H)$ when $H$ is separable.

The spatial formulas in \eqref{intro:standard} are common to these results, but the domains impose different coherence requirements. Each operator-norm classification uses one implementation throughout its prescribed fixed-rank stratum or infinite-rank, infinite-corank sector. Likewise, an isometry between individual $L^p$-perturbation classes has one implementation on its entire source class. The extended $L^p$-metric assigns infinite distance to distinct classes, so their implementations may be chosen independently whenever the classes can be matched isometrically.

\section{Common algebra and metric linearization}\label{sec:common}
All inner products are linear in the first variable. For $x,y\in H$, the operator $x\otimes y$ is defined by $z\mapsto\langle z,y\rangle x$. For an idempotent $E$, write $R(E)$ for the orthogonal projection onto $\Ran E$.

\begin{lemma}[Block form and range projection]\label{common:block}
Relative to $H=\Ran E\oplus(\Ran E)^\perp$, every idempotent has the form
\[
E=\begin{pmatrix}I&A\\0&0\end{pmatrix}.
\]
The operator $E+E^*-I$ is invertible and
\begin{equation}\label{common:rangeformula}
R(E)=E(E+E^*-I)^{-1}.
\end{equation}
If $E-F\in L^p(\B(H),\tr)$, then $R(E)-R(F)\in L^p(\B(H),\tr)$.
\end{lemma}
\begin{proof}
The block form follows from $E|_{\Ran E}=I$. In this decomposition,
\[
E+E^*-I=\begin{pmatrix}I&A\\A^*&-I\end{pmatrix},\qquad
(E+E^*-I)^2=\begin{pmatrix}I+AA^*&0\\0&I+A^*A\end{pmatrix}.
\]
The square is bounded below by $I$, and block multiplication gives \eqref{common:rangeformula}. The final assertion follows from the ideal property of $L^p(\B(H),\tr)$ and the resolvent identity for the two inverses.
\end{proof}

\begin{lemma}[Corners, charts, and square-zero lines]\label{common:corners}
Let $\mathcal A$ be either $\B(H)$ with its operator norm or $L^p(\B(H),\tr)$ with its $L^p$-norm. The tangent space at $E$ is
\[
X_E=\{X\in\mathcal A:EX+XE=X\}=V_E\oplus W_E,
\]
where $V_E=E\mathcal A(I-E)$ and $W_E=(I-E)\mathcal A E$. The canonical projection onto $V_E$ is $Q_EX=EX(I-E)$. Local similarity coordinates are
\begin{equation}\label{common:charts}
\Theta_E(X)=e^{[X,E]}Ee^{-[X,E]},\qquad D\Theta_E(0)X=X.
\end{equation}
An affine line $E+tN$ consists of idempotents for every real $t$ if and only if
\begin{equation}\label{common:squarezero}
EN+NE=N,\qquad N^2=0.
\end{equation}
These conditions also imply that $E+zN$ is an idempotent for every $z\in\C$.
\end{lemma}
\begin{proof}
Differentiating $E(t)^2=E(t)$ gives the tangent equation. It implies $EXE=0$ and $(I-E)X(I-E)=0$, and hence the corner decomposition. The identity $[[X,E],E]=X$ proves the derivative assertion in \eqref{common:charts}. The analytic inverse function theorem on the split idempotent manifold gives local coordinates; see \cite{ref3}. Similarity preserves rank and corank. When $E$ has finite rank, all the displayed tangent operators have rank at most $2\rank E$, so the same local construction works for every $L^p$-norm. Finally,
\[
(E+tN)^2-(E+tN)=t(EN+NE-N)+t^2N^2,
\]
which proves the last two assertions.
\end{proof}

\begin{lemma}[Metric blow-up]\label{common:blowup}
Let $\varphi$ be a surjective isometry between two idempotent manifolds with the local charts above, and set $F=\varphi(E)$. For a free ultrafilter $\mathcal U$ on $\N$, the real Banach ultrapowers $(X_E)_{\mathcal U}$ and $(X_F)_{\mathcal U}$ admit a surjective real-linear isometry.
\end{lemma}
\begin{proof}
Choose $t_n=n^{-2}$. For all sufficiently large $n$, the local charts define
\[
f_n(X)=t_n^{-1}\Theta_F^{-1}\bigl(\varphi(\Theta_E(t_nX))\bigr)
\]
on the ball $\|X\|\leq n$, since $t_nn\to0$. Shrink the chart neighborhoods so that both charts satisfy
\[
\tfrac12\|X-Y\|\leq\|\Theta(X)-\Theta(Y)\|\leq\tfrac32\|X-Y\|.
\]
In particular, $\|f_n(X)\|\leq3\|X\|$ wherever the construction is used sufficiently close to the chart origins. The chart expansions and the isometry property give, uniformly on every fixed bounded ball,
\[
\left|\|f_n(X)-f_n(Y)\|-\|X-Y\|\right|\longrightarrow0,
\qquad f_n(0)=0.
\]
Here all norms are taken in the chosen ambient space. The domains have radii tending to infinity, and changing finitely many coordinates has no effect in the ultrapower. Thus $[X_n]\mapsto[f_n(X_n)]$ is a well-defined isometry fixing zero.

For surjectivity, apply the same construction to $\varphi^{-1}$ and any bounded sequence $(Y_n)$ in $X_F$. The inverse sequence is bounded by the same local first-order estimates. The forward and inverse coordinate formulas compose exactly to the identity for all sufficiently large indices; in particular, $f_n(X_n)-Y_n\to0$. The resulting ultrapower isometry is onto, and the real Mazur--Ulam theorem \cite{ref11} makes it real linear. See \cite{ref8} for the standard construction of Banach ultrapowers.
\end{proof}

\begin{remark}\label{common:ultrapowerremark}
The blow-up lemma gives linear isometries between ultrapowers; it does not assert differentiability of $\varphi$. In finite rank, a further identification with concrete corners on the Hilbert ultrapower is available. The operator-norm argument below instead uses internal rank-one elements in both rank regimes; the finite-rank $L^p$ argument uses the concrete identification.
\end{remark}

\section{\texorpdfstring{The $L^\infty$ metric: tangent geometry and projection preservation}{The Linfinity metric: tangent geometry and projection preservation}}
\label{sec:operator-projections}

In this section, $\mathcal I$ denotes either $\Id_k(H)$, where $1\leq k<\infty$ and $H$ is an arbitrary infinite-dimensional complex Hilbert space, or $\Id_{\infty,\infty}(H)$, where $H$ is separable. Write
\[
 \mathcal P=\{p\in\mathcal I:p=p^*\}.
\]
We use the $L^\infty$ norm on $\B(H)$, that is, the operator norm, throughout this section and Section~\ref{sec:operator-classification}. Inner products are linear in the first variable, and $(u\otimes v)\xi=\langle\xi,v\rangle u$.

In each case the domain is the specified stratum or sector; these sections do not classify maps defined on the whole set $\Id(H)$.

\begin{theorem}[$L^\infty$ projection preservation]
\label{op:projection-preservation}
Every surjective $L^\infty$-isometry $\varphi:\mathcal I\to\mathcal I$ satisfies $\varphi(\mathcal P)=\mathcal P$.
\end{theorem}

\subsection{Similarity coordinates and rank-one directions}
Fix $e\in\mathcal I$ and put
\[
 X_e=\{x\in\B(H):ex+xe=x\},\qquad
 V_e=e\B(H)(I-e),\qquad W_e=(I-e)\B(H)e.
\]
The common local-coordinate construction identifies $X_e$ with the tangent space at $e$; its analytic local similarity charts have derivative the identity. We shall use the decomposition
\begin{equation}
 X_e=V_e\oplus W_e,\qquad Q_ex=ex(I-e).
 \label{op:corner-decomposition}
\end{equation}
Here the ambient space in that construction is $L^\infty(\B(H),\tr)=\B(H)$.

Recall that a real-linear projection $P$ on a real Banach space $X$ is an $M$-projection if
\begin{equation}
 \norm{x}=\max\{\norm{Px},\norm{(I-P)x}\}\qquad(x\in X).
 \label{op:M-projection}
\end{equation}
It is nontrivial when $P\ne0,I$; see~\cite{ref7}.

\begin{lemma}[Rank-one dichotomy]
\label{op:rank-one-dichotomy}
Let $X$ be a real-linear subspace of $\B(K)$ for a complex Hilbert space $K$, and let $P$ be an $M$-projection on $X$. If $x\in X$ is nonzero and has rank one, then $Px=0$ or $Px=x$.
\end{lemma}
\begin{proof}
Put $a=Px$, $b=(I-P)x$, $\alpha=\norm{a}$, and $\beta=\norm{b}$. Then
\begin{equation}
 \norm{sa+tb}=\max\{|s|\alpha,|t|\beta\}\qquad(s,t\in\R).
 \label{op:rank-one-max}
\end{equation}
Suppose that $\alpha,\beta>0$. If $\alpha>\beta$, then $\norm{x}=\alpha$. Write $x=\alpha u\otimes v$ with unit vectors $u,v$. The point $\lambda=1$ is interior to the interval on which $\norm{a+\lambda b}=\alpha$, so $\norm{x+tb}=\alpha$ for all sufficiently small real $t$. Applying this inequality to $v$, squaring, and adding the inequalities for $t$ and $-t$ gives $bv=0$. Applying the same argument to $x^*,b^*$ gives $b^*u=0$. This argument uses only the norm identity; membership of the adjoints in $X$ is unnecessary. Relative to the orthogonal domain and range decompositions determined by $v,u$, respectively, we obtain
\[
 \norm{x+tb}=\max\{\alpha,|t|\beta\}\qquad(t\in\R).
\]
Take $\lambda=-\alpha/\beta$ and $t=\lambda-1$. Since $x+tb=a+\lambda b$, equation~\eqref{op:rank-one-max} gives $\alpha$, whereas the last formula gives $\alpha+\beta$, a contradiction. The case $\beta>\alpha$ is symmetric.

If $\alpha=\beta$, write $x=\alpha u\otimes v$ with unit vectors $u,v$ and put $y=a-b$. Then
\[
 \norm{x}=\norm{y}=\alpha,\qquad \norm{x+y}=\norm{x-y}=2\alpha.
\]
For unit vectors $\xi_n$ with $\norm{(x+y)\xi_n}\to2\alpha$, the triangle inequality forces $\norm{x\xi_n}\to\alpha$. After changing phases, $\xi_n\to v$. Equality in the Hilbert-space triangle inequality then gives $yv=\alpha u$. Applying the same reasoning to $x-y$ gives $-yv=\alpha u$, a contradiction. Thus $\alpha\beta=0$.
\end{proof}

\subsection{An ultrapower characterization of orthogonal projections}
Fix a free ultrafilter $\mathcal U$ on $\mathbb N$. For a Banach space $X$, let
\[
 X_{\mathcal U}=\ell^\infty(X)\big/\{(x_n):\lim_{\mathcal U}\norm{x_n}=0\},
 \qquad \norm{[(x_n)]}=\lim_{\mathcal U}\norm{x_n}.
\]
We use the usual Hilbert space ultrapower $H_{\mathcal U}$ and the standard ultrapower construction~\cite{ref8}. The map
\[
 J_e:(X_e)_{\mathcal U}\longrightarrow\B(H_{\mathcal U}),
 \qquad J_e([(x_n)])[(\xi_n)]=[(x_n\xi_n)]
\]
is an isometric embedding. The upper norm estimate is immediate; for the reverse estimate choose unit $\xi_n$ with $\norm{x_n\xi_n}\geq\norm{x_n}-1/n$. We identify $(X_e)_{\mathcal U}$ with this image, without asserting that the image is the full tangent space on $H_{\mathcal U}$.

Let $\mathcal R_e^+$ consist of the classes $[(z_n)]$ with $z_n\in V_e$, $\operatorname{rank}z_n=1$, and $\lim_{\mathcal U}\norm{z_n}>0$. Define $\mathcal R_e^-$ similarly using $W_e$. These will be called nonzero internal rank-one elements.

\begin{lemma}
\label{op:internal-connectivity}
Each element of $\mathcal R_e^+\cup\mathcal R_e^-$ has rank one as an operator on $H_{\mathcal U}$. Each of the sets $\mathcal R_e^+$ and $\mathcal R_e^-$ is path connected.
\end{lemma}
\begin{proof}
Write $R=eH$ and $N=\ker e$. Every rank-one element of $V_e$ can be written
\[
 z_n=\rho_nu_n\otimes a_n,\qquad \rho_n=\norm{z_n}>0,
\]
where $u_n\in R$, $a_n\in R^\perp$, and $\norm{u_n}=\norm{a_n}=1$. If $\rho=\lim_{\mathcal U}\rho_n>0$, then, for every bounded sequence $(\xi_n)$,
\[
 [(z_n\xi_n)]=\left(\lim_{\mathcal U}\rho_n\langle\xi_n,a_n\rangle\right)[(u_n)].
\]
The image is one-dimensional and is nonzero, as seen by using $[(a_n)]$. For $W_e$ use $N,N^\perp$.

We record a uniform path construction valid in every nonzero complex Hilbert space $L$. Given unit vectors $u,v$, choose $\theta\in[-\pi,\pi]$ such that $\langle u,e^{i\theta}v\rangle\geq0$ is real. The normalized line segment from $u$ to $e^{i\theta}v$ has denominator at least $1/\sqrt2$, and hence a Lipschitz constant bounded by an absolute constant. Follow it by the phase path from $e^{i\theta}v$ to $v$. After reparametrizing onto $[0,1]$, all these paths have a common Lipschitz bound. This includes $\dim L=1$.

First interpolate $\rho_n$ linearly to $1$. This joins $[(z_n)]$ to $[(u_n\otimes a_n)]$ through nonzero internal rank-one elements, since $(1-t)\rho+t>0$. Fix unit vectors $u_0\in R$ and $a_0\in R^\perp$. Apply the uniform paths just constructed to join $u_n$ to $u_0$ and $a_n$ to $a_0$. Their tensor products define a continuous path in the ultrapower because
\[
 \norm{u\otimes a-u'\otimes a'}\leq\norm{u-u'}+\norm{a-a'}
\]
for unit vectors. This joins the class to $u_0\otimes a_0$. The same proof with $N,N^\perp$ treats $\mathcal R_e^-$. All four spaces used here are nonzero because $e$ has nonzero range and nonzero kernel.
\end{proof}

\begin{lemma}[Failure of the corner maximum]
\label{op:corner-failure}
If $e\ne e^*$, there are rank-one $v\in V_e$ and $w\in W_e$ such that
\[
 \norm{v-w}>\max\{\norm{v},\norm{w}\}.
\]
\end{lemma}
\begin{proof}
Let $P_R,P_N$ denote the orthogonal projections onto $R=eH$ and $N=\ker e$, and put $T=P_R|_N$. Since $e$ is not orthogonal, $c=\norm{T}>0$. For $n\in N$,
\[
 \operatorname{dist}(n,R)\geq\frac{\norm{n}}{\norm{I-e}},
\]
because $(I-e)(n-r)=n$ for $r\in R$. Thus $c<1$. Take unit approximate eigenvectors $n_j\in N$ for the top spectral value $c^2$ of $T^*T$, and put
\[
 a_j=\norm{Tn_j},\qquad r_j=Tn_j/a_j,\qquad d_j=\norm{P_Nr_j}.
\]
Then $a_j,d_j\to c$, since $T^*=P_N|_R$; also $\langle r_j,n_j\rangle=a_j$. For large $j$ define
\[
 \alpha_j=\frac{P_{R^\perp}n_j}{\norm{P_{R^\perp}n_j}^2},\quad
 \beta_j=\frac{P_{N^\perp}r_j}{\norm{P_{N^\perp}r_j}^2},\quad
 v_j=r_j\otimes\alpha_j,\quad w_j=n_j\otimes\beta_j.
\]
Then $v_j\in V_e$, $w_j\in W_e$, and
\[
 \norm{v_j}=(1-a_j^2)^{-1/2},\qquad
 \norm{w_j}=(1-d_j^2)^{-1/2}.
\]
Furthermore, $(v_j-w_j)(r_j-n_j)=-(r_j+n_j)$, so
\[
 \norm{v_j-w_j}\geq\sqrt{\frac{1+a_j}{1-a_j}}
 \longrightarrow\frac{1+c}{\sqrt{1-c^2}},
\]
whereas $\max\{\norm{v_j},\norm{w_j}\}\to(1-c^2)^{-1/2}$. Choose one sufficiently large $j$.
\end{proof}

\begin{lemma}[Internal rank-one perturbation]
\label{op:internal-perturbation}
For every nonzero $x\in(X_e)_{\mathcal U}$ there is $z\in\mathcal R_e^+\cup\mathcal R_e^-$ such that
\[
 \norm{z}<\norm{x}<\norm{x+z}.
\]
\end{lemma}
\begin{proof}
Represent $x=[(x_n)]$ with $\norm{x_n}\leq M$, and put $m=\norm{x}>0$. Choose unit vectors $\xi_n$ with $\norm{x_n\xi_n}\geq\norm{x_n}-1/n$. Set
\[
 x_n^+=ex_n(I-e),\quad x_n^-=(I-e)x_ne,\quad
 y_n=x_n\xi_n,\quad y_n^\pm=x_n^\pm\xi_n.
\]
Since $y_n=y_n^++y_n^-$, at least one of
\[
 \operatorname{Re}\langle y_n^+,y_n\rangle\geq\tfrac12\norm{y_n}^2,
 \qquad
 \operatorname{Re}\langle y_n^-,y_n\rangle\geq\tfrac12\norm{y_n}^2
\]
holds. One choice holds on a set in $\mathcal U$. Suppose it is the plus choice, and intersect that set with $\{n:\norm{y_n}>m/2\}$ to obtain $A\in\mathcal U$. For $n\in A$ put $s_n=(I-e)\xi_n\ne0$ and
\[
 a_n=\frac{P_{R^\perp}s_n}{\norm{P_{R^\perp}s_n}^2},\qquad
 z_n^0=y_n^+\otimes a_n.
\]
Then $z_n^0\in V_e$, $z_n^0\xi_n=y_n^+$, and
\[
 \norm{z_n^0}=\frac{\norm{y_n^+}}{\operatorname{dist}(s_n,R)}
 \leq\norm{e}\norm{I-e}M=:C,
\]
using $\operatorname{dist}(s_n,R)\geq\norm{s_n}/\norm{I-e}$ and $\norm{y_n^+}\leq\norm{e}M\norm{s_n}$. Also
\[
 \norm{z_n^0}\geq\norm{y_n^+}\geq\norm{y_n}/2>m/4.
\]
Choose $t>0$ with $tC<m$, set $z_n=tz_n^0$ on $A$, and use $z_n=q/n$ off $A$, for a fixed nonzero rank-one $q\in V_e$. Thus $z=[(z_n)]\in\mathcal R_e^+$ and $0<tm/4\leq\norm{z}<m$. On $A$,
\[
 \norm{(x_n+z_n)\xi_n}^2
 =\norm{y_n+ty_n^+}^2\geq(1+t)\norm{y_n}^2.
\]
Taking the ultralimit gives $\norm{x+z}\geq\sqrt{1+t}\,m>m$.

For the minus choice, put $r_n=e\xi_n$ and replace $a_n$ by $P_{N^\perp}r_n/\norm{P_{N^\perp}r_n}^2$ and $y_n^+$ by $y_n^-$. The estimates
\[
 \operatorname{dist}(r_n,N)\geq\norm{r_n}/\norm{e},\qquad
 \norm{y_n^-}\leq\norm{I-e}M\norm{r_n}
\]
give the same bound and yield $z\in\mathcal R_e^-$.
\end{proof}

\begin{proposition}[Intrinsic tangent invariant]
\label{op:tangent-invariant}
For every $e\in\mathcal I$ and every free ultrafilter $\mathcal U$,
\[
 e=e^*\quad\Longleftrightarrow\quad
 (X_e)_{\mathcal U}\text{ admits a nontrivial real $M$-projection}.
\]
\end{proposition}
\begin{proof}
If $e=p=p^*$, then, relative to $pH\oplus p^\perp H$,
\[
 x=\begin{pmatrix}0&B\\C&0\end{pmatrix}\in X_p,
 \qquad \norm{x}=\max\{\norm{B},\norm{C}\}.
\]
Hence $Q_p$ is a nontrivial $M$-projection and induces one on the ultrapower.

Conversely, let $e\ne e^*$ and suppose $P$ is a nontrivial $M$-projection on $(X_e)_{\mathcal U}$. By Lemmas~\ref{op:rank-one-dichotomy} and~\ref{op:internal-connectivity}, $P$ either fixes or annihilates each internal rank-one element, and path connectedness makes its choice constant on each of $\mathcal R_e^+$ and $\mathcal R_e^-$. If $P$ fixes both sets, take $0\ne b\in\ker P$ and use Lemma~\ref{op:internal-perturbation} to obtain an internal rank-one $z$ with
\[
 \norm{b+z}>\max\{\norm{b},\norm{z}\},
\]
contrary to the $M$-identity. If $P$ annihilates both sets, use instead $0\ne a\in\operatorname{Ran}P$. In either mixed case, regard the $v,w$ from Lemma~\ref{op:corner-failure} as constant ultrapower elements. The $M$-identity would give $\norm{v-w}=\max\{\norm{v},\norm{w}\}$, again a contradiction.
\end{proof}

\begin{proof}[Proof of Theorem~\ref{op:projection-preservation}]
Let $e=e^*$ and $f=\varphi(e)$. Proposition~\ref{op:tangent-invariant} provides a nontrivial $M$-projection $P$ on $(X_e)_{\mathcal U}$. Conjugating $P$ by the real-linear isometry furnished by the common metric linearization lemma gives a nontrivial $M$-projection on $(X_f)_{\mathcal U}$. Proposition~\ref{op:tangent-invariant} implies $f=f^*$. Applying the same argument to $\varphi^{-1}$ proves equality of the projection sets.
\end{proof}

\section{\texorpdfstring{The $L^\infty$ metric: extension from projections and classification}{The Linfinity metric: extension from projections and classification}}
\label{sec:operator-classification}

The common extension argument applies to either of the two spaces $\mathcal I$ specified in Section~\ref{sec:operator-projections}. For an idempotent $E$, write $R(E)$ for the orthogonal projection onto its range, and define
\[
 \mathcal D=\{E\in\mathcal I:\operatorname{rank}(E-R(E))=1\}.
\]
Thus every member of $\mathcal D$ has exactly one nonorthogonal two-dimensional block.

\subsection{Two perturbation formulas}
We shall use the following standard isolated-eigenvalue formulas~\cite{ref9}, with their short proofs included to specify signs and the inner-product convention.

\begin{lemma}
\label{op:eigenvalue-variation}
Let $A(t)$ be a $C^2$ curve of self-adjoint operators. If $\lambda$ is an isolated simple eigenvalue of $A(0)$ with unit eigenvector $x$, and $\lambda(t)$ is its local continuation, then
\[
 \lambda''(0)=\langle A''(0)x,x\rangle
 +2\left\langle(\lambda I-A(0))^{-1}_{x^\perp}A'(0)x,A'(0)x\right\rangle.
\]
The reduced inverse is extended by zero on $\C x$.
\end{lemma}
\begin{proof}
Choose a $C^2$ unit eigenvector branch $x(t)$ with $x(0)=x$ and $\langle x'(0),x\rangle=0$. Differentiating the eigenvalue equation gives
\[
 \lambda'(0)=\langle A'(0)x,x\rangle,\qquad
 x'(0)=(\lambda I-A(0))^{-1}_{x^\perp}A'(0)x.
\]
Differentiating a second time and taking the inner product with $x$ gives the asserted formula. The second term is real since the reduced inverse is self-adjoint.
\end{proof}

\begin{lemma}
\label{op:rank-one-asymptotic}
For $X\in\B(H_1,H_2)$, unit vectors $\alpha\in H_2$, $\beta\in H_1$, and $|\lambda|=1$,
\[
 \norm{X-t\lambda\alpha\otimes\beta}
 =t-\operatorname{Re}\bigl(\overline\lambda\langle X\beta,\alpha\rangle\bigr)
 +O(t^{-1})\qquad(t\to\infty).
\]
\end{lemma}
\begin{proof}
Put $R=\alpha\otimes\beta$ and $\varepsilon=t^{-1}$. Then
\[
 t^{-2}(X-t\lambda R)^*(X-t\lambda R)
 =\beta\otimes\beta-\varepsilon(\lambda X^*R+\overline\lambda R^*X)
 +\varepsilon^2X^*X.
\]
The simple isolated largest eigenvalue $1$ at $\varepsilon=0$ continues as the largest eigenvalue
\[
 \mu(\varepsilon)=1-2\varepsilon\operatorname{Re}
 \bigl(\overline\lambda\langle X\beta,\alpha\rangle\bigr)+O(\varepsilon^2).
\]
Taking $t\sqrt{\mu(t^{-1})}$ proves the formula.
\end{proof}

\subsection{Rigidity on one nonorthogonal block}
\begin{lemma}
\label{op:one-block-rigidity}
Let $\psi:\mathcal I\to\mathcal I$ be a surjective $L^\infty$-isometry fixing every element of $\mathcal P$. Then, for each $E\in\mathcal D$, $\psi(E)=E$ or $\psi(E)=E^*$.
\end{lemma}
\begin{proof}
Write $p=R(E)$. In the orthogonal decomposition $pH\oplus p^\perp H$,
\[
 E=\begin{pmatrix}I&A\\0&0\end{pmatrix},\qquad A=au\otimes\xi,
\]
where $a>0$ and $u\in pH$, $\xi\in p^\perp H$ are unit vectors. Let $S=\operatorname{span}\{u,\xi\}$ and $N=S^\perp$. Then
\[
 E=E_a\oplus P,\qquad E_a=\begin{pmatrix}1&a\\0&0\end{pmatrix},
\]
where $P$ is an orthogonal projection on $N$ of rank $k-1$ in the finite-rank case, and of infinite rank and corank in the other case. Every $Q\oplus P$, $Q\in\Proj_1(S)$, belongs to $\mathcal P$.

Set $s=\sqrt{1+a^2}$,
\[
 c=\sqrt{\frac{s+1}{2s}},\qquad d=\sqrt{\frac{s-1}{2s}},\qquad
 w_0=cu+d\xi,\quad w_1=-du+c\xi,\quad Q_j=w_j\otimes w_j.
\]
In the basis $w_0,w_1$,
\[
 E_a=\begin{pmatrix}1+\alpha&\beta\\-\beta&-\alpha\end{pmatrix},
 \qquad \alpha=(s-1)/2,\quad\beta=a/2.
\]
Put $\rho=\alpha+\beta>0$ and $q_j=Q_j\oplus P$. Direct computation of the two singular values gives
\begin{equation}
 \norm{E-q_0}=\rho,\qquad\norm{E-q_1}=1+\rho.
 \label{op:extreme-projections}
\end{equation}
For $F=\psi(E)$, $T=F-q_0$, and $J=q_0-q_1$, these equalities read
\[
 \norm{T}=\rho,\qquad \norm{T+J}=1+\rho,\qquad \norm{J}=1.
\]
Choose unit $x_n$ with $\norm{(T+J)x_n}\to1+\rho$. The triangle inequality implies $\norm{Jx_n}\to1$. Since $J^*J=P_S$ and $S$ is two-dimensional, a subsequence converges in norm to a unit vector $x\in S$. Equality in the triangle inequality then gives $Tx=\rho Jx$ and $T^*Tx=\rho^2x$.

Put $u_0=Q_0x$ and $v_0=Q_1x$. We have
\[
 Tx=\rho(u_0-v_0),\qquad T^*(u_0-v_0)=\rho(u_0+v_0).
\]
The idempotency of $F=q_0+T$ gives $T^2+q_0T+Tq_0=T$. Applying it to $x$ and using the preceding identities yields
\[
 Tu_0=A_0u_0-B_0v_0,\qquad Tv_0=B_0u_0-A_0v_0,
 \qquad A_0=\frac{\rho^2}{1+2\rho},\quad
 B_0=\frac{\rho(1+\rho)}{1+2\rho}.
\]
Taking the inner product of $T^*(u_0-v_0)=\rho(u_0+v_0)$ with $u_0$ shows
\[
 A_0\norm{u_0}^2+B_0\norm{v_0}^2=\rho\norm{u_0}^2.
\]
Since $\rho-A_0=B_0$, we obtain $\norm{u_0}=\norm{v_0}=1/\sqrt2$. Thus these two vectors span $S$ and $F(S)\subset S$. Consequently
\[
 F=\begin{pmatrix}F_0&Z\\0&D\end{pmatrix}\quad\text{on }S\oplus N.
\]
In the orthonormal basis $\sqrt2u_0,\sqrt2v_0$,
\[
 F_0=\begin{pmatrix}1+A_0&B_0\\-B_0&-A_0\end{pmatrix}.
\]
Using $s^2-a^2=1$ gives $A_0=\alpha$ and $B_0=\beta$. In particular $F_0$ is unitarily similar to $E_a$.

For every $Q\in\Proj_1(S)$,
\[
 \norm{F_0-Q}\leq\norm{F-(Q\oplus P)}=\norm{E_a-Q}.
\]
Integrate over $\Proj_1(S)$ with its unitarily invariant probability measure, which has full support. Unitary similarity of $F_0,E_a$ makes the two integrals equal, so continuity makes the displayed inequality an equality for every $Q$. In the fixed basis $w_0,w_1$, the previous formulas give, for a real $\theta$,
\[
 F_0=\begin{pmatrix}1+\alpha&\beta e^{-i\theta}\\-\beta e^{i\theta}&-\alpha\end{pmatrix}.
\]
For $Q_+=\tfrac12\begin{pmatrix}1&1\\1&1\end{pmatrix}$, the matrices $F_0-Q_+$ and $E_a-Q_+$ have the same $L^2$ norm for the matrix trace, and
\[
 |\det(E_a-Q_+)|^2=\tfrac14,\qquad
 |\det(F_0-Q_+)|^2=\tfrac14+\beta^2\sin^2\theta.
\]
Their equal $L^\infty$ norms determine the larger singular value; their equal $L^2$ norms then determine the smaller one. Their determinant moduli agree, so $\sin\theta=0$. Therefore
\begin{equation}
 F_0=E_a\quad\text{or}\quad F_0=E_a^*.
 \label{op:two-block-choices}
\end{equation}

For a unit left singular vector $y_Q$ corresponding to the largest singular value of $F_0-Q$, apply $(F-(Q\oplus P))^*$ to $(y_Q,0)$. The equality $\norm{F-(Q\oplus P)}=\norm{F_0-Q}$ gives $Z^*y_Q=0$. These vectors span $S$: take $Q=R(F_0)$ and $Q=I_S-R(F_0)$. In coordinates in which $F_0=\begin{pmatrix}1&a\\0&0\end{pmatrix}$, the first left singular vector is the first coordinate vector, whereas that vector is not an eigenvector of
\[
 \begin{pmatrix}1&a\\0&-1\end{pmatrix}
 \begin{pmatrix}1&0\\a&-1\end{pmatrix}
 =\begin{pmatrix}1+a^2&-a\\-a&1\end{pmatrix}.
\]
Thus $Z=0$ and $F=F_0\oplus D$.

Let $G=F_0\oplus P$ and $C=D-P$. By~\eqref{op:two-block-choices}, $G=E$ or $E^*$, and, for every $q\in\mathcal P$,
\begin{equation}
 \norm{F-q}=\norm{G-q}.
 \label{op:distance-identity}
\end{equation}
Using $q_0$ in this identity gives $\norm{C}\leq\rho<a$, since
\[
 2(a-\rho)=a+1-\sqrt{1+a^2}>0.
\]
Let $p_0=R(F_0)$ and choose orthogonal unit $\gamma,\eta\in S$ so that $F_0-p_0=a\gamma\otimes\eta$. Put $\bar q=p_0\oplus P$. Then
\[
 F-\bar q=a\gamma\otimes\eta\oplus C.
\]
The top eigenvalue $a^2$ of $(F-\bar q)^*(F-\bar q)$ is isolated and simple, with unit eigenvector $\eta$.

For a unit $z\in P^\perp N$, replace the projection direction $\gamma$ by $\gamma(t)=\cos t\,\gamma+\sin t\,z$:
\[
 q_z(t)=\bar q-\gamma\otimes\gamma+\gamma(t)\otimes\gamma(t).
\]
This is a unitary conjugate of $\bar q$, hence belongs to $\mathcal P$. For $H_C(t)=(F-q_z(t))^*(F-q_z(t))$, differentiation gives
\[
 H_C'(0)\eta=-az,\qquad \langle H_C''(0)\eta,\eta\rangle=0.
\]
By the spectral gap and Lemma~\ref{op:eigenvalue-variation}, the second derivative of $\norm{F-q_z(t)}^2$ at zero is
\[
 2a^2\langle(a^2I-C^*C)^{-1}z,z\rangle.
\]
For $G$ it is $2$. Equality of the distance functions gives
\[
 a^2\langle(a^2I-C^*C)^{-1}z,z\rangle=1.
\]
Expanding the inverse in a norm-convergent series, all terms beyond the constant term are nonnegative. In particular $\norm{Cz}^2=0$; therefore $CP^\perp=0$.

For completeness, take a unit $\ell\in PN$, when this space is nonzero, and replace $\ell$ by $\ell(t)=\cos t\,\ell+\sin t\,\eta$:
\[
 q_\ell(t)=\bar q-\ell\otimes\ell+\ell(t)\otimes\ell(t).
\]
Now $H_C'(0)\eta=-C^*\ell$ and $\langle H_C''(0)\eta,\eta\rangle=2$. Equality with the derivative for $G$ yields
\[
 2+2\langle(a^2I-C^*C)^{-1}C^*\ell,C^*\ell\rangle=2.
\]
Strict positivity of the inverse gives $C^*P=0$; if $PN=0$ this identity is automatic. Finally $F^*,G^*$ have the same distances to every $q\in\mathcal P$. The first rotation argument, applied to their two-dimensional blocks and exterior error $C^*$, gives $C^*P^\perp=0$. Hence $C^*=0$ and $D=P$. Together with~\eqref{op:two-block-choices}, this proves $F=E$ or $F=E^*$.
\end{proof}

\begin{lemma}
\label{op:constant-choice}
Under the assumptions of Lemma~\ref{op:one-block-rigidity}, the choice is constant on $\mathcal D$: either $\psi(E)=E$ for all $E\in\mathcal D$, or $\psi(E)=E^*$ for all $E\in\mathcal D$.
\end{lemma}
\begin{proof}
The two alternatives give disjoint relatively open subsets of $\mathcal D$. Indeed, if $\psi(E)=E$ and $\psi(F)=F^*$, then
\[
 \norm{E-E^*}\leq\norm{E-F^*}+\norm{F^*-E^*}=2\norm{E-F}.
\]
The reverse choice is handled in the same way.

The set $\mathcal D$ is path connected. Every member is uniquely $p+A$, with $p\in\mathcal P$ and $A\in p\B(H)(I-p)$ nonzero of rank one. For $p,q\in\mathcal P$, choose a unitary $V$ with $VpV^*=q$. In finite rank it may be chosen on the finite-dimensional space $pH+qH$ and extended by the identity. In infinite rank its existence follows by choosing unitaries between the respective ranges and kernels. A bounded self-adjoint logarithm of $V$ provides a norm-continuous unitary path $V_t$ from $I$ to $V$. Thus $V_t(p+A)V_t^*$ joins $p+A$ to $q+VAV^*$. The nonzero rank-one operators in $q\B(H)(I-q)$ form a path-connected set, by writing them as $ru\otimes v$, $r>0$, with unit vectors in $qH$ and $q^\perp H$. This supplies the rest of the path. Connectedness now forces one of the two relatively open subsets to be empty.
\end{proof}

\subsection{Recovering an arbitrary idempotent}
\begin{lemma}
\label{op:arbitrary-recovery}
If $\psi$ fixes every element of $\mathcal P\cup\mathcal D$, then $\psi(E)=E$ for all $E\in\mathcal I$.
\end{lemma}
\begin{proof}
Fix $E\in\mathcal I$ and write, relative to $pH\oplus p^\perp H$,
\[
 p=R(E),\qquad E=\begin{pmatrix}I&A\\0&0\end{pmatrix},\qquad F=\psi(E).
\]
For unit $\alpha\in pH$, $\xi\in p^\perp H$, and $c\in\C$, put
\[
 G_c=\begin{pmatrix}I&0\\c\xi\otimes\alpha&0\end{pmatrix}.
\]
This is an idempotent with the same rank and corank as $p$. For $c\ne0$, its sole nonorthogonal block on $\operatorname{span}\{\alpha,\xi\}$ is $\begin{pmatrix}1&0\\c&0\end{pmatrix}$; subtracting its range projection gives a nonzero matrix of rank one. Thus $G_c\in\mathcal D$ and $\psi(G_c)=G_c$. The off-diagonal blocks of $E-G_c$ have orthogonal domains and ranges, so
\[
 \norm{E-G_c}=\max\{\norm{A},|c|\}.
\]
For $|c|>\norm{A}$ this yields $\norm{F\alpha-\alpha-c\xi}\leq|c|$. Writing $F\alpha-\alpha=\lambda\xi+\eta$, $\eta\perp\xi$, gives
\[
 |\lambda-c|^2+\norm{\eta}^2\leq|c|^2
\]
for all sufficiently large complex $c$. Choosing $c$ on the ray opposite to $\lambda$ first forces $\lambda=0$, then $\eta=0$. Hence $F|_{pH}=I$, and
\[
 F=\begin{pmatrix}I&B\\0&D\end{pmatrix}.
\]
We eliminate $D$ by an argument valid in both rank regimes. For $|c|>\norm{A}$, the operator $T=F-G_c$ satisfies $T\alpha=-c\xi$ and $\norm{T}=|c|$. Positivity of $\norm{T}^2I-T^*T$ implies $T^*T\alpha=|c|^2\alpha$. The $p^\perp H$ component of this identity is $-cD^*\xi=0$. Since $\xi$ was arbitrary, $D=0$.

For unit $\beta\in p^\perp H$ put
\[
 K_c=\begin{pmatrix}I&c\alpha\otimes\beta\\0&0\end{pmatrix}.
\]
Again $K_c\in\mathcal D$ for $c\ne0$, and isometry gives
\[
 \norm{A-c\alpha\otimes\beta}=\norm{B-c\alpha\otimes\beta}\qquad(c\in\C).
\]
Apply Lemma~\ref{op:rank-one-asymptotic} with $c=t\lambda$ and let $t\to\infty$. Taking $\lambda=1,i$ recovers the real and imaginary parts of $\langle A\beta,\alpha\rangle$ and $\langle B\beta,\alpha\rangle$. They agree for all $\alpha,\beta$, so $A=B$ and $F=E$.
\end{proof}

\begin{proposition}
\label{op:projection-extension}
A surjective $L^\infty$-isometry $\psi:\mathcal I\to\mathcal I$ fixing every orthogonal projection is either $\psi(E)=E$ for every $E$, or $\psi(E)=E^*$ for every $E$.
\end{proposition}
\begin{proof}
Lemmas~\ref{op:one-block-rigidity} and~\ref{op:constant-choice} give one of the two choices throughout $\mathcal D$. If it is the identity, apply Lemma~\ref{op:arbitrary-recovery}. Otherwise apply that lemma to $E\mapsto\psi(E)^*$, which is again a surjective $L^\infty$-isometry fixing $\mathcal P$ and now fixes $\mathcal D$.
\end{proof}

\subsection{\texorpdfstring{The two $L^\infty$ classifications}{The two Linfinity classifications}}
\begin{theorem}[Finite-rank $L^\infty$ classification]
\label{op:finite-rank-classification}
Let $H$ be any infinite-dimensional complex Hilbert space and $1\leq k<\infty$. For every surjective $L^\infty$-isometry $\varphi:\Id_k(H)\to\Id_k(H)$, there exists a unitary or antiunitary $U$ such that one of the following formulas holds throughout $\Id_k(H)$:
\[
 \varphi(E)=UEU^*,\qquad\text{or}\qquad\varphi(E)=UE^*U^*.
\]
Conversely, each displayed map is a surjective $L^\infty$-isometry.
\end{theorem}
\begin{proof}
By Theorem~\ref{op:projection-preservation}, the restriction to $\Proj_k(H)$ is a surjective Grassmann isometry. Its classification~\cite{ref1,ref5} gives $\varphi(p)=UpU^*$ for a unitary or antiunitary $U$. The complementary-projection alternative is excluded because $I-p$ has infinite rank. Normalize by $\psi(E)=U^*\varphi(E)U$ and apply Proposition~\ref{op:projection-extension}. The converses follow from invariance of the $L^\infty$ norm under these operations.
\end{proof}

\begin{theorem}[Infinite rank and infinite corank]
\label{op:infinite-rank-classification}
Let $H$ be separable and infinite-dimensional. Every surjective $L^\infty$-isometry of $\Id_{\infty,\infty}(H)$ has one of the following four forms, with one unitary or antiunitary $U$ and one formula used globally:
\[
 \begin{aligned}
 \varphi(E)&=UEU^*,&\quad \varphi(E)&=UE^*U^*,\\
 \varphi(E)&=U(I-E)U^*,&\quad \varphi(E)&=U(I-E^*)U^*.
 \end{aligned}
\]
Each of these forms defines a surjective $L^\infty$-isometry.
\end{theorem}
\begin{proof}
Projection preservation and the infinite Grassmann classification~\cite{ref6} give either $\varphi(p)=UpU^*$ for every $p\in\Proj_{\infty,\infty}(H)$, or $\varphi(p)=U(I-p)U^*$ for every such $p$. In the first case normalize by $U^*\varphi(\cdot)U$; in the second normalize by $I-U^*\varphi(\cdot)U$. Complementation preserves $\Id_{\infty,\infty}(H)$, since it exchanges the range and kernel, and is an $L^\infty$-isometry. The normalized map fixes all projections, so Proposition~\ref{op:projection-extension} supplies the two possibilities in each case. The converse is immediate.
\end{proof}

\section{\texorpdfstring{$L^p$}{Lp}-isometries on idempotents of fixed finite rank}
\label{sec:fixed-rank}

Throughout this section, $H$ is an infinite-dimensional complex Hilbert space and
$k\geq 1$ is finite. We work in the noncommutative $L^p$-space
$L^p(\mathcal B(H),\operatorname{Tr})$ associated with the canonical trace on
$\mathcal B(H)$, with
\[
 \|A\|_p=\bigl(\operatorname{Tr}|A|^p\bigr)^{1/p},
 \qquad 1\leq p<\infty.
\]
The objects under consideration are
\[
 \mathcal I_k(H)=\{E\in\mathcal B(H):E^2=E,\ \operatorname{rank}E=k\},
 \qquad
 \mathcal P_k(H)=\{P\in\mathcal I_k(H):P=P^*\}.
\]
Every difference of two elements of $\mathcal I_k(H)$ has rank at most $2k$.
Consequently,
\[
 d_p(E,F)=\|E-F\|_p
 \qquad(E,F\in\mathcal I_k(H))
\]
is a finite-valued metric. We call an isometry for this metric an
$L^p$-isometry. For a finite-rank operator $A$, let $R(A)$ denote the
orthogonal projection onto its range, with $R(0)=0$.

\begin{theorem}\label{fr:classification}
Let $1\leq p<\infty$. Every surjective $L^p$-isometry
$\varphi:\mathcal I_k(H)\to\mathcal I_k(H)$ has one of the forms
\[
 \varphi(E)=UEU^*\quad(E\in\mathcal I_k(H)),
 \qquad\text{or}\qquad
 \varphi(E)=UE^*U^*\quad(E\in\mathcal I_k(H)),
\]
where $U$ is unitary or antiunitary. No differentiability assumption is
needed.
\end{theorem}

The proof separates into three ranges. For $p\ne2$, an intrinsic splitting
property of the tangent spaces identifies the orthogonal projections. In the
strictly convex range $1<p<\infty$, the normalized isometry then extends along
affine families of idempotents. At the $L^1$ endpoint we instead use nearest
projections, while the $L^2$ case is treated by extending the isometry to the
ambient real Hilbert space.

\subsection{Projection preservation for \texorpdfstring{$p\ne2$}{p different from 2}}

We use the following disjointness criterion in
$L^p(\mathcal B(H),\operatorname{Tr})$: if $p\ne2$ and
\begin{equation}\label{fr:disjointness}
 \|A+B\|_p^p=\|A-B\|_p^p=\|A\|_p^p+\|B\|_p^p,
\end{equation}
then $A^*B=AB^*=0$. For $1<p<\infty$, this is the equality case of the
Clarkson--McCarthy inequalities \cite{ref12,ref2}. At $p=1$, it follows from
the equality case of the $L^1$ triangle inequality \cite{ref10}. Indeed,
equality for $A+B$ gives $A=VP$ and $B=VQ$, where $P,Q\geq0$ and $V$ is a
common compatible polar partial isometry. Equality for $A-B$ then gives
\[
 \|P-Q\|_1=\operatorname{Tr}P+\operatorname{Tr}Q,
\]
which forces the supports of $P$ and $Q$ to be orthogonal: the self-adjoint
contraction $\operatorname{sgn}(P-Q)$ must be $I$ on the support of $P$ and
$-I$ on the support of $Q$. Thus $PQ=0$, proving the criterion.

Recall from Lemma~\ref{common:corners} that the tangent space at $E\in\mathcal I_k(H)$ is
\[
 X_E=\{X\in L^p(\mathcal B(H),\operatorname{Tr}):EX+XE=X\}
     =V_E\oplus W_E,
\]
where
\[
 V_E=E L^p(\mathcal B(H),\operatorname{Tr})(I-E),\qquad
 W_E=(I-E)L^p(\mathcal B(H),\operatorname{Tr})E,
 \qquad Q_EX=EX(I-E).
\]
All operators in $X_E$ have rank at most $2k$.

A real-linear projection $T:Y\to Y$ on a real normed subspace of
$L^p(\mathcal B(H),\operatorname{Tr})$ is called a \emph{$p$-splitting
projection} if
\[
 \|y\|_p^p=\|Ty\|_p^p+\|(I-T)y\|_p^p\qquad(y\in Y).
\]
This identity implies that $T$ is bounded.

\begin{lemma}\label{fr:splitting}
Let $1\leq p<\infty$ and $p\ne2$. Every $p$-splitting projection on $X_E$
belongs to $\{0,I,Q_E,I-Q_E\}$. Moreover, $Q_E$ is $p$-splitting if and only
if $E=E^*$. Consequently,
\begin{equation}\label{fr:projection-criterion}
 E=E^*
 \quad\Longleftrightarrow\quad
 (X_E,\|\cdot\|_p)\text{ admits a nontrivial $p$-splitting projection}.
\end{equation}
\end{lemma}

\begin{proof}
Let $T$ be a $p$-splitting projection and let $X\in X_E$ have rank one. Put
$A=TX$ and $B=(I-T)X$. Applying the splitting identity to $A+B$ and $A-B$
gives \eqref{fr:disjointness}. Their left supports are orthogonal, as are
their right supports, so
\[
 \operatorname{rank}(A+B)=\operatorname{rank}A+\operatorname{rank}B.
\]
It follows that $TX=X$ or $TX=0$.

Relative to the algebraic decomposition $H=EH\dotplus\ker E$, every
element of $X_E$ has the off-diagonal form
\[
 \begin{pmatrix}0&B\\ C&0\end{pmatrix},
\]
of rank $\operatorname{rank}B+\operatorname{rank}C$. Hence the nonzero
rank-one vectors of $X_E$ form two classes, lying respectively in $V_E$ and
$W_E$. Each class is path connected and spans its corner over $\mathbb R$.
The alternatives $TX=X$ and $TX=0$ define disjoint relatively closed and
open subsets of each class. They are therefore constant on each class,
which gives the four stated possibilities.

If $E=P=P^*$, the two off-diagonal corners have orthogonal left and right
supports, so their $p$th norm powers add. Conversely, suppose that $EH$ is
not orthogonal to $\ker E$. Choose $r\in EH$ and $n\in\ker E$ with
$\langle n,r\rangle\ne0$, together with nonzero vectors
\[
 a\in(EH)^\perp=\ker E^*,\qquad
 b\in(\ker E)^\perp=\operatorname{Ran}E^*.
\]
Then $v=r\otimes a\in V_E$ and $w=n\otimes b\in W_E$, whereas
\[
 v^*w=\langle n,r\rangle a\otimes b\ne0.
\]
A splitting identity for $Q_E$, applied to $v+w$ and $v-w$, would contradict
\eqref{fr:disjointness}.
\end{proof}

The metric blow-up of Lemma~\ref{common:blowup} transports this intrinsic criterion without
requiring a derivative of $\varphi$. Write
\[
 \widehat L_E:(X_E)_{\mathcal U}\longrightarrow
                   (X_{\varphi(E)})_{\mathcal U}
\]
for the surjective real-linear isometry supplied by that construction.
Finite rank allows us to realize these Banach ultrapowers as concrete
tangent spaces on the Hilbert ultrapower.

\begin{lemma}\label{fr:ultrapower}
Let $H_{\mathcal U}$ be the Hilbert ultrapower and let $E_{\mathcal U}$ be
the operator induced by the constant sequence $E$. Then
$\operatorname{rank}E_{\mathcal U}=k$, and
\[
 (X_E)_{\mathcal U}\cong X_{E_{\mathcal U}}
      \subset L^p(\mathcal B(H_{\mathcal U}),\operatorname{Tr})
\]
isometrically as real normed spaces.
\end{lemma}

\begin{proof}
The range of $E$ is finite dimensional, so its bounded vector sequences have
ultralimits in that same range. Thus $E_{\mathcal U}$ has rank $k$. For a
bounded sequence $X_n\in X_E$, define
\[
 \Psi([(X_n)])[(\xi_n)]=[(X_n\xi_n)].
\]
The tangent equation passes to the ultrapower. Since
$\operatorname{rank}X_n\leq2k$, the min--max principle for singular values
gives
\[
 s_j\bigl(\Psi([(X_n)])\bigr)
      =\lim_{\mathcal U}s_j(X_n)\qquad(1\leq j\leq2k).
\]
One can also obtain this by choosing singular value decompositions
termwise and taking ultrapowers of the finitely many singular vectors.
Consequently,
\[
 \|\Psi([(X_n)])\|_p=\lim_{\mathcal U}\|X_n\|_p.
\]

For surjectivity, take $X\in X_{E_{\mathcal U}}$. Its two corners have rank
at most $k$, so $X$ is a finite sum of rank-one operators on
$H_{\mathcal U}$. Choose bounded vector representatives to obtain a bounded
sequence $Y_n$ of finite-rank operators inducing $X$, and set
\[
 X_n=EY_n(I-E)+(I-E)Y_nE.
\]
Then $X_n\in X_E$, their ranks and $L^p$-norms are uniformly bounded, and
they induce $X$. This corner construction also accounts for the dual
vectors when $E$ is not self-adjoint.
\end{proof}

\begin{proposition}\label{fr:projection-preservation}
For $p\ne2$, every surjective $L^p$-isometry of $\mathcal I_k(H)$ maps
$\mathcal P_k(H)$ onto itself.
\end{proposition}

\begin{proof}
For $P\in\mathcal P_k(H)$, transport a nontrivial splitting projection on
$X_{P_{\mathcal U}}$ through Lemma~\ref{common:blowup} and Lemma~\ref{fr:ultrapower}.
Lemma~\ref{fr:splitting}, applied on $H_{\mathcal U}$, shows that
$\varphi(P)_{\mathcal U}$ is self-adjoint. A constant operator is
self-adjoint on its ultrapower exactly when it is self-adjoint on $H$.
Apply the same argument to $\varphi^{-1}$.
\end{proof}

By the Grassmann-space $L^p$-isometry theorem
\cite[Theorem~1.2]{ref13}, the restriction to $\mathcal P_k(H)$ is
implemented by a unitary or antiunitary. The complementary alternative
that can occur in a finite factor does not arise here, because $H$ is
infinite dimensional and $k$ is finite. After spatial normalization, we may
therefore assume that
\begin{equation}\label{fr:normalization}
 \varphi(P)=P\qquad(P\in\mathcal P_k(H)).
\end{equation}

\subsection{Extension along affine wings in the strictly convex range}

Assume $1<p<\infty$ and $p\ne2$, and fix $P\in\mathcal P_k(H)$.
By Lemma~\ref{fr:splitting}, the blow-up $\widehat L_P$ either preserves
$(V_P)_{\mathcal U}$ and $(W_P)_{\mathcal U}$ or interchanges them. For
$v\in V_P$, the curve $\Theta_P(t(v+v^*))$ consists of orthogonal
projections, because $[v+v^*,P]$ is skew-adjoint. Equation
\eqref{fr:normalization} therefore gives
\[
 \widehat L_P(v+v^*)=v+v^*,
\]
where the vectors are identified with their constant ultrapower classes.
Uniqueness of the corner decomposition implies, for every $v\in V_P$,
\[
 \widehat L_P(v)=v\quad\text{in the preserving case},\qquad
 \widehat L_P(v)=v^*\quad\text{in the interchanging case}.
\]

For $0\leq t\leq1$, the points $P+tv$ lie in $\mathcal I_k(H)$ and satisfy
equality in the triangle inequality between $P$ and $P+v$. Strict convexity
of $L^p(\mathcal B(H),\operatorname{Tr})$ gives
\[
 \varphi(P+tv)=P+t\bigl(\varphi(P+v)-P\bigr).
\]
The blow-up of this curve at zero has value $\varphi(P+v)-P$. Hence one
fixed alternative holds on the entire row wing at $P$:
\[
 \varphi(P+v)=P+v\quad(v\in V_P),\qquad\text{or}\qquad
 \varphi(P+v)=P+v^*\quad(v\in V_P).
\]
The column wing has the same alternative, by the corresponding action of
$\widehat L_P$.

We now make the alternative independent of the base projection; the same
argument will also be used at $p=1$. Suppose $P,Q\in\mathcal P_k(H)$ and
$\|P-Q\|<1$, where $\|\cdot\|$ is the operator norm. Then
\[
 A=Q|_{PH}:PH\longrightarrow QH
\]
is invertible, and
\[
 E=A^{-1}Q
\]
is an idempotent with $R(E)=P$ and $R(E^*)=Q$. Thus $E$ belongs to the row
wing at $P$ and the column wing at $Q$. If $P\ne Q$, then $E\ne E^*$, so
the alternatives at $P$ and $Q$ must agree. The Grassmann space
$\mathcal P_k(H)$ is path connected. Subdividing a path into pieces whose
adjacent projections have operator-norm distance less than one makes the
alternative independent of $P$. Every idempotent belongs to the row wing
at its range projection. This proves Theorem~\ref{fr:classification} for
$1<p<\infty$, $p\ne2$.

\subsection{The \texorpdfstring{$L^1$}{L1} endpoint and nearest projections}

Assume $p=1$ and retain \eqref{fr:normalization}. Strict convexity is
unavailable, so we recover an idempotent from its distances to the rank-$k$
orthogonal projections. Write
\[
 E=\begin{pmatrix}I&A\\0&0\end{pmatrix}
\]
relative to $R(E)H\oplus R(E)^\perp H$. For every $S\in\mathcal P_k(H)$,
\[
 2\|E-S\|_1
   =\|E-S\|_1+\|S-E^*\|_1
   \geq\|E-E^*\|_1=2\|A\|_1,
\]
and equality is attained at $S=R(E)$. Denote the set of minimizers by
$M(E)$.

Put $r=\operatorname{rank}A$. A singular value decomposition gives an
orthogonal decomposition
\[
 H=M\oplus H_1\oplus H_2\oplus N,\qquad
 E=I_M\oplus E_D\oplus0_N,\qquad
 E_D=\begin{pmatrix}I&D\\0&0\end{pmatrix},
\]
where $\dim M=k-r$, $\dim H_1=\dim H_2=r$, and $D>0$ after identifying
$H_1$ with $H_2$. For self-adjoint $X$, set
\[
 C_X=(I+X^2)^{-1},\qquad
 G(X)=\begin{pmatrix}C_X&C_XX\\XC_X&X^2C_X\end{pmatrix}.
\]

\begin{proposition}\label{fr:nearest-projections}
For this decomposition,
\[
 M(E)=\{I_M\oplus G(X)\oplus0_N:0\leq X\leq D\}.
\]
\end{proposition}

\begin{proof}
Suppose $S\in M(E)$. Equality in the $L^1$ triangle inequality gives common
compatible positive polar factors for $E-S$ and $S-E^*$. On the active
block,
\[
 E_D-E_D^*=JK_0,\qquad
 J=\begin{pmatrix}0&I\\-I&0\end{pmatrix},\qquad
 K_0=\operatorname{diag}(D,D).
\]
The positive factors sum to $K_0$ on the active block and zero on
$M\oplus N$. Thus they vanish on $M\oplus N$ and on all cross blocks.
Uniqueness of the polar partial isometry on the active block gives
\[
 S=I_M\oplus S_0\oplus0_N,\qquad
 E_D-S_0=JB,\qquad 0\leq B\leq K_0.
\]
Write
\[
 S_0=\begin{pmatrix}C&Z\\Z^*&F\end{pmatrix}.
\]
Then
\[
 B=-J(E_D-S_0)
   =\begin{pmatrix}Z^*&F\\I-C&D-Z\end{pmatrix}.
\]
Positivity gives $Z=Z^*\geq0$ and $F=I-C$. Since $S_0$ is a projection,
\[
 C^2+Z^2=C,\qquad CZ=ZC.
\]
If $Cu=0$, then $Zu=0$; positivity of $B$ forces $(I-C)u=0$, hence $u=0$.
In finite dimension $C$ is therefore invertible. Define $X=C^{-1}Z\geq0$.
The projection equation gives $C=(I+X^2)^{-1}$ and $S_0=G(X)$. Now
\[
 B=\begin{pmatrix}XC&X^2C\\X^2C&D-CX\end{pmatrix},\qquad
 B\begin{pmatrix}-Xu\\u\end{pmatrix}
   =\begin{pmatrix}0\\(D-X)u\end{pmatrix}.
\]
Thus $D-X\geq0$.

Conversely, suppose $0\leq X\leq D$, and put $C=(I+X^2)^{-1}$. For
$S_0=G(X)$, the same matrix $B=-J(E_D-S_0)$ satisfies
\[
 B=\begin{pmatrix}I\\X\end{pmatrix}
      XC\begin{pmatrix}I&X\end{pmatrix}
      +\begin{pmatrix}0&0\\0&D-X\end{pmatrix}\geq0.
\]
Since $J$ is unitary,
\[
 \|E_D-S_0\|_1=\operatorname{Tr}B=\operatorname{Tr}D.
\]
This is the lower bound already established, proving the assertion.
\end{proof}

For a subset $\mathcal A$ of a normed space, use the Bouligand tangent cone
\[
 \operatorname{Tan}(\mathcal A,a)
 =\{Z:t_n\downarrow0,\ Z_n\to Z,\ a+t_nZ_n\in\mathcal A
                      \text{ for some }(t_n),(Z_n)\}.
\]
A cone $C$ is \emph{pointed} if $C\cap(-C)=\{0\}$.

\begin{lemma}\label{fr:pointed-cones}
For a finite-dimensional positive definite $D$, the tangent cone of $[0,D]$
is pointed exactly at $0$ and $D$.
\end{lemma}

\begin{proof}
At the endpoints the cones consist respectively of positive and negative
semidefinite operators. Congruence by $D^{-1/2}$ reduces the other points
to the interval $[0,I]$. If $Y$ has an eigenvalue strictly between $0$ and
$1$, the corresponding rank-one projection and its negative are tangent
directions. Otherwise $Y\ne0,I$ is a nontrivial projection. Rotating a
range vector toward a kernel vector gives a nonzero tangent direction and
its negative. The cone is therefore not pointed.
\end{proof}

The graph map $X\mapsto G(X)$ is a real-analytic parametrization, with
inverse $X=S_{21}S_{11}^{-1}$ on its image. Its derivative identifies the
corresponding tangent cones. By Proposition~\ref{fr:nearest-projections}
and Lemma~\ref{fr:pointed-cones}, for non-self-adjoint $E$ the only points
of $M(E)$ with pointed tangent cones are
\[
 R(E)=I_M\oplus G(0)\oplus0_N,\qquad
 R(E^*)=I_M\oplus G(D)\oplus0_N.
\]
Consequently, if
\[
 \|E-S\|_1=\|F-S\|_1\qquad(S\in\mathcal P_k(H)),
\]
then $M(E)=M(F)$ and
\[
 \{R(E),R(E^*)\}=\{R(F),R(F^*)\}.
\]
If one idempotent is an orthogonal projection, taking that projection as
$S$ gives equality of the idempotents directly. Otherwise, the ordered
range projections determine the range and kernel, since
$\ker E=R(E^*)^\perp H$. Thus in all cases $F=E$ or $F=E^*$. Together with
\eqref{fr:normalization}, this gives
\begin{equation}\label{fr:l1-pointwise}
 \varphi(E)\in\{E,E^*\}\qquad(E\in\mathcal I_k(H)).
\end{equation}

Fix $P\in\mathcal P_k(H)$. If nonzero $A,B\in V_P$ have opposite choices
in \eqref{fr:l1-pointwise}, preservation of distances gives
\[
 \|A-B\|_1=\|A-B^*\|_1=\|A\|_1+\|B\|_1.
\]
For fixed nonzero $A$, this is impossible when $B$ is sufficiently close
to $A$. The choice is therefore locally constant on $V_P\setminus\{0\}$,
which is path connected, and hence is constant on the row wing.
Injectivity forces $E$ and $E^*$ to have the same choice, so the column
wing at $P$ has that choice as well. The nearby-range-projection
construction in the preceding subsection makes the choice independent of
$P$. This proves Theorem~\ref{fr:classification} at $p=1$.

\subsection{The \texorpdfstring{$L^2$}{L2} case}

For $p=2$, every tangent space is a real Hilbert space, so the splitting
criterion cannot distinguish orthogonal projections. Instead, we extend
the isometry to the ambient real Hilbert space
$L^2(\mathcal B(H),\operatorname{Tr})$, with inner product
$\operatorname{Re}\operatorname{Tr}(A^*B)$.

Let
\[
 \mathcal F_0(H)=\{A\in\mathcal B(H):\operatorname{rank}A<\infty,
                                      \ \operatorname{Tr}A=0\}.
\]
We first observe that
\begin{equation}\label{fr:real-span}
 \operatorname{span}_{\mathbb R}
       (\mathcal I_k(H)-\mathcal I_k(H))=\mathcal F_0(H),
 \qquad
 \overline{\mathcal F_0(H)}^{\|\cdot\|_2}
       =L^2(\mathcal B(H),\operatorname{Tr}).
\end{equation}
Indeed, a rank-one square-zero operator $N=x\otimes y$ is the difference
$(P+N)-P$ for a rank-$k$ projection $P$ with $x\in PH$ and $y\perp PH$.
Such nilpotents span $\mathcal F_0(H)$ over $\mathbb R$: the off-diagonal
matrix units and their complex multiples are available, while on every
two-dimensional coordinate space
\[
 \begin{pmatrix}1&0\\0&-1\end{pmatrix}
 =\frac12\begin{pmatrix}1&1\\-1&-1\end{pmatrix}
  +\frac12\begin{pmatrix}1&-1\\1&-1\end{pmatrix}
\]
is a sum of rank-one square-zero matrices. Conversely, differences of
rank-$k$ idempotents have trace zero. For density, if $A$ has finite rank
and trace $t$, choose $n$ orthogonal rank-one projections $q_j$ outside its
finite-dimensional support. Then
\[
 A_n=A-\frac{t}{n}\sum_{j=1}^nq_j\in\mathcal F_0(H),\qquad
 \|A_n-A\|_2=\frac{|t|}{\sqrt n}.
\]

Fix $E_0\in\mathcal I_k(H)$. Polarization of the preserved distances shows
that
\[
 E-E_0\longmapsto\varphi(E)-\varphi(E_0)
\]
preserves all real pairwise inner products. It therefore extends linearly
to a well-defined isometry on their real span. By \eqref{fr:real-span}, it
has a unique continuous extension
\[
 L:L^2(\mathcal B(H),\operatorname{Tr})
       \longrightarrow L^2(\mathcal B(H),\operatorname{Tr}).
\]
Surjectivity of $\varphi$ makes the closed range of $L$ dense. Thus $L$ is
surjective, and
\begin{equation}\label{fr:l2-affine}
 \varphi(E)=L(E)+b\qquad(E\in\mathcal I_k(H))
\end{equation}
for a fixed $b\in L^2(\mathcal B(H),\operatorname{Tr})$.

For $0\leq q\leq k$, set
\[
 \mathcal C_q=\{N\in L^2(\mathcal B(H),\operatorname{Tr}):
                      N^2=0,\ \operatorname{rank}N\leq q\}.
\]
An operator $N$ belongs to $\mathcal C_k$ exactly when some affine line
$E+\mathbb RN$ is contained in $\mathcal I_k(H)$. For the forward direction,
choose a $k$-dimensional space between $\operatorname{Ran}N$ and $\ker N$
and let $E$ be its orthogonal projection. For the converse, comparison of
coefficients in $(E+tN)^2=E+tN$ gives
\[
 N^2=0,\qquad EN+NE=N.
\]
In the algebraic off-diagonal form
$N=\begin{pmatrix}0&B\\C&0\end{pmatrix}$, we have
$\operatorname{Ran}B\subseteq\ker C\subseteq EH$. Consequently,
\[
 \operatorname{rank}N=\operatorname{rank}B+\operatorname{rank}C\leq k.
\]
Equation \eqref{fr:l2-affine} and its inverse imply
\begin{equation}\label{fr:preserve-square-zero-cone}
 L(\mathcal C_k)=\mathcal C_k.
\end{equation}

\begin{lemma}\label{fr:square-zero-tangent}
Let $1\leq q\leq k$ and $N\in\mathcal C_q$. The cone
$\operatorname{Tan}(\mathcal C_q,N)$ is a real-linear space if and only if
$\operatorname{rank}N=q$.
\end{lemma}

\begin{proof}
Suppose $\operatorname{rank}N=q$. Put
\[
 R=\operatorname{Ran}N,\qquad D=(\ker N)^\perp,\qquad
 M=\ker N\ominus R.
\]
In the decomposition $R\oplus M\oplus D$, the only nonzero block of $N$
is its invertible block $N_{13}$. A fixed bounded similarity, differing
from the identity by finite rank, normalizes this block to $I$ after
identifying $D$ with $R$. Similarity preserves $\mathcal C_q$ and linearity
of its tangent cones, so assume this normalization.

If $Z$ is tangent, differentiation of the square-zero equation gives
$NZ+ZN=0$. We also have $Z(\ker N)\subseteq\operatorname{Ran}N$. Indeed,
otherwise choose $u\in\ker N$ with $Zu\notin\operatorname{Ran}N$.
Compress $N+t_nZ_n$ to the columns $D\oplus\mathbb Cu$ and the rows $R$
together with the nonzero component of $Zu$ in $R^\perp$. After dividing
the last column by $t_n$, these $(q+1)$-dimensional compressions converge
to an invertible matrix, contradicting their rank bound. The two necessary
conditions give exactly
\[
 Z=\begin{pmatrix}
      Z_{11}&Z_{12}&Z_{13}\\
      0&0&Z_{23}\\
      0&0&-Z_{11}
    \end{pmatrix}.
\]
Every such $Z$ equals $[X,N]$: take
\[
 X_{31}=-Z_{11},\qquad X_{32}=-Z_{12},\qquad
 X_{21}=Z_{23},\qquad X_{11}-X_{33}=Z_{13}.
\]
The curve $e^{tX}Ne^{-tX}$ remains in $\mathcal C_q$ and has derivative
$Z$. The tangent cone is therefore the displayed real-linear space.

Now let $r=\operatorname{rank}N<q$ and $s=q-r$. Outside the
finite-dimensional support of $N$, choose mutually orthogonal spaces
$D_1,R_1,D_2,R_2$, each of dimension $s$. Let $X$ and $Y$ be isomorphisms
$D_1\to R_1$ and $D_2\to R_2$, extended by zero. Then
$N+tX,N+tY\in\mathcal C_q$, so $X$ and $Y$ are tangent directions. If
$X+Y$ were tangent, compress to the columns $D_0\oplus D_1\oplus D_2$ and
the rows $\operatorname{Ran}N\oplus R_1\oplus R_2$, where
$N|_{D_0}$ is an isomorphism onto its range. After division of the last two
column groups by $t_n$, these matrices converge to an invertible matrix of
size
\[
 r+2s=2q-r>q.
\]
This contradiction proves that the tangent cone is not linear.
\end{proof}

Since $L$ preserves tangent cones, \eqref{fr:preserve-square-zero-cone} and
Lemma~\ref{fr:square-zero-tangent} first identify the rank-$k$ stratum and
then its complement $\mathcal C_{k-1}$. Induction gives
$L(\mathcal C_q)=\mathcal C_q$ for every $q\leq k$. In particular, $L$
preserves rank-one square-zero operators in both directions. We record the
linear-preserver consequence that will be used next.

\begin{lemma}\label{fr:linear-preserver}
Let $T:\mathcal F_0(H)\to\mathcal F_0(H)$ be a bijective real-linear map,
with $T$ and $T^{-1}$ bounded for the $L^2$-norm. If $T$ preserves rank-one
square-zero operators in both directions, then
\[
 T(A)=cSAS^{-1}\quad(A\in\mathcal F_0(H)),
 \qquad\text{or}\qquad
 T(A)=cSA^*S^{-1}\quad(A\in\mathcal F_0(H)),
\]
where $c\ne0$ and $S$ is a bounded invertible linear or conjugate-linear
operator on $H$.
\end{lemma}

\begin{proof}
Let $\mathcal N_1$ denote the nonzero rank-one square-zero operators. For
$N,M\in\mathcal N_1$, real linearity gives
\[
 N+M\in\mathcal N_1
 \quad\Longleftrightarrow\quad
 T(N)+T(M)\in\mathcal N_1.
\]
Apply the rank-one nilpotent relation-preserver theorem of Du, Hou and Bai
\cite{ref4} to this restriction. Its projective conclusion gives a
semilinear similarity or a dual semilinear similarity, with a field
automorphism and, initially, nonzero scalar factors depending on $N$.

Continuity excludes discontinuous field automorphisms. Indeed, for
orthogonal $e_1,e_2,e_3$, apply $T$ to $e_1\otimes(e_2+ze_3)$ and take a
ratio of two projective coordinates in the fixed semilinear
representation. The point-dependent factor cancels, and continuity of
$T$ implies continuity of the field automorphism. It is therefore either
the identity or complex conjugation.

The resulting projective map is continuous. For example, when range lines
map to range lines, if $x_n\to x\ne0$, choose $y\perp x$, $y\ne0$, and
project $y$ onto $x_n^\perp$ to obtain $y_n\to y$. Then
\[
 T(x_n\otimes y_n)\longrightarrow T(x\otimes y)\ne0
\]
gives convergence of the range lines. In the dual case the same argument
uses kernel functionals, or adjoints. The inverse projective map is
continuous as well.

Let $S$ implement the resulting linear or conjugate-linear projective map.
To see that $S$ is bounded, suppose $x_n\to x$ and $Sx_n\to y$, and put
$w=y-Sx$. For every $v\ne0$,
\[
 [v+x_n-x]\longrightarrow[v].
\]
Projective continuity implies either $Sv+w=0$ or $[Sv+w]=[Sv]$; hence
$w\in\mathbb CSv$. Two linearly independent images force $w=0$. The
closed graph and bounded inverse theorems now apply. After the Riesz
identification in the dual case, the fixed model is
\[
 F(A)=SAS^{-1}\qquad\text{or}\qquad F(A)=SA^*S^{-1}.
\]

Write $T(N)=\lambda_NF(N)$. If $N,M,N+M\in\mathcal N_1$ and $N,M$ are
complex-linearly independent, additivity gives
\[
 \lambda_N=\lambda_M=\lambda_{N+M}.
\]
These relations connect all of $\mathcal N_1$. For $N=x\otimes u$ and
$M=y\otimes v$, choose $h\perp x,y$, avoiding the finitely many prohibited
lines, and use
\[
 x\otimes u,\qquad x\otimes h,\qquad
 y\otimes h,\qquad y\otimes v.
\]
Adjacent independent terms have the same range or the same kernel
functional, and their sum is rank one and square zero. Proportional terms
are bridged by an independent term with the same range. These choices are
possible because $H$ is infinite dimensional. Thus $\lambda_N=c$ is
constant. Equation \eqref{fr:real-span} and real linearity extend the
formula to $\mathcal F_0(H)$.
\end{proof}

By \eqref{fr:real-span}, $L$ maps $\mathcal F_0(H)$ bijectively onto
itself, so Lemma~\ref{fr:linear-preserver} applies. Its $L^2$-isometry
property forces $S$ to be a scalar multiple of a unitary or antiunitary.
In the linear similarity case, for orthogonal unit vectors $x,y$,
\[
 1=|c|\,\|Sx\|\,\|(S^{-1})^*y\|.
\]
For fixed $y$, the value of $\|Sx\|$ is constant on the unit sphere of
$y^\perp$. Any two unit vectors have a common orthogonal unit vector, so
$\|Sx\|=\alpha\|x\|$ for every $x$. The other cases reduce to this one by
adjunction or a fixed conjugation. Cancelling $\alpha>0$ gives
\[
 L(A)=cJ(A)\qquad(A\in\mathcal F_0(H)),\qquad |c|=1,
\]
where $J(A)=UAU^*$ or $J(A)=UA^*U^*$.

Since $E-E_0\in\mathcal F_0(H)$, the affine formula becomes
\[
 \varphi(E)=cJ(E)+D,\qquad D=\varphi(E_0)-cJ(E_0),
\]
with $D$ finite rank. As $P=J(E)$ ranges over $\mathcal I_k(H)$, the
operator $cP+D$ is idempotent. Subtracting the equations for $P,Q$ and
using \eqref{fr:real-span} yields
\[
 AD+DA=(1-c)A\qquad(A\in\mathcal F_0(H)).
\]
Thus $C=D-\frac{1-c}{2}I$ anticommutes with every off-diagonal matrix
unit. Comparing matrix entries first makes every off-diagonal entry of
$C$ zero; for distinct $i,j$, its diagonal entries satisfy
$c_{ii}+c_{jj}=0$. Three distinct indices force every diagonal entry to
vanish. Hence
\[
 D=\frac{1-c}{2}I.
\]
Since $D$ has finite rank and $H$ is infinite dimensional, $c=1$ and $D=0$.
This proves Theorem~\ref{fr:classification} for $p=2$.

\section{\texorpdfstring{$L^p$ perturbation classes and the extended metric}{Lp perturbation classes and the extended metric}}\label{sec:classes}

Throughout this section, $H$ is a separable infinite-dimensional complex Hilbert space and $1\leq p<\infty$. We work in the canonical type-I noncommutative $L^p$ space
\[
 L^p=L^p(\mathcal B(H),\operatorname{Tr}),
 \qquad \|A\|_p=\bigl(\operatorname{Tr}(|A|^p)\bigr)^{1/p},
\]
where $\operatorname{Tr}$ is the canonical trace. Put
\[
 \mathcal I(H)=\{E\in\mathcal B(H):E^2=E\},
 \qquad \mathcal C_p(E_0)=\{E\in\mathcal I(H):E-E_0\in L^p\}.
\]
The extended $L^p$ distance on $\mathcal I(H)$ is
\[
 d_p(E,F)=
 \begin{cases}
 \|E-F\|_p,&E-F\in L^p,\\
 +\infty,&E-F\notin L^p.
 \end{cases}
\]
Thus the $L^p$ perturbation classes $\mathcal C_p(E_0)$ are the finite-distance equivalence classes of this metric. They need not be connected. The principal rigidity statement concerns a whole finite-distance class, including all of its connected components.

\begin{theorem}[Rigidity of an $L^p$ perturbation class]\label{cl:main}
Suppose that $E_0,F_0\in\mathcal I(H)$ have infinite rank and infinite corank. Every surjective $L^p$ isometry
\[
 \Phi:\mathcal C_p(E_0)\longrightarrow\mathcal C_p(F_0)
\]
has one of the forms
\begin{equation}\label{cl:standard-forms}
 \Phi(E)=UEU^*,\quad UE^*U^*,\quad U(I-E)U^*,\quad U(I-E^*)U^*,
\end{equation}
where $U$ is one fixed unitary or antiunitary operator. Conversely, any map in \eqref{cl:standard-forms} that sends $E_0$ into $\mathcal C_p(F_0)$ restricts to a surjective $L^p$ isometry between these classes.
\end{theorem}

The proof recovers the relative Fredholm index from the metric and then normalizes rank-one orthogonal increments. The local argument for $1<p<\infty$ uses strict convexity of $L^p$, whereas the $L^1$ endpoint requires a separate norm-equality argument. We describe the entire extended metric space after establishing this classwise result.

\subsection{Finite-dimensional blocks and the relative index}

We first record a finite-rank modification that remains available when a perturbation class contains no orthogonal projection. If $A,A+q$ are idempotents and $q$ is an orthogonal projection, expansion gives $Aq+qA=0$. Comparing blocks for $qH\oplus q^\perp H$ then gives $Aq=qA=0$.

\begin{lemma}[Installation of a prescribed orthogonal block]\label{cl:install}
Let $E_0$ have infinite rank and infinite corank. For any finite-dimensional subspace $V\subset H$ and any orthogonal projection $p_V$ on $V$, there is an invertible $S=I+K$, with $K$ of finite rank, such that $G=SE_0S^{-1}$ is orthogonally reduced by $V$ and $G|_V=p_V$. Moreover, $G-E_0$ has finite rank, and an $L^p$-continuous path of idempotents joins $E_0$ to $G$.
\end{lemma}

\begin{proof}
Write $R=\operatorname{Ran}E_0$, $N=\ker E_0$, $r=\operatorname{rank}p_V$, and $s=\dim V-r$. Choose subspaces $R_0\subset R$ and $N_0\subset N$ of dimensions $r,s$, together with bounded projections $T_R:R\to R_0$ and $T_N:N\to N_0$. Then
\[
 T=T_RE_0+T_N(I-E_0)
\]
is a finite-rank idempotent commuting with $E_0$, of rank $\dim V$. Choose a finite-dimensional $L$ containing $V$, $\operatorname{Ran}T$, and $\operatorname{Ran}T^*$. It reduces $T$ orthogonally, and $T|_{L^\perp}=0$. The decompositions
\[
 L=R_0\dotplus N_0\dotplus\ker(T|_L),\qquad
 L=p_VV\oplus(I_V-p_V)V\oplus(L\ominus V)
\]
have corresponding summands of equal dimensions. Choose an invertible $S_L$ taking the first decomposition to the second, respecting its summands, and extend it by the identity on $L^\perp$. Then $STS^{-1}=P_V$, the orthogonal projection onto $V$. Since this projection commutes with $SE_0S^{-1}$, $V$ reduces $G$; the first two summands give the prescribed values of $G|_V$. The difference $SE_0S^{-1}-E_0$ is finite rank. Finally, the complex finite-dimensional general linear group is path connected. A path from $I_L$ to $S_L$, extended by the identity, gives the asserted $L^p$-continuous path.
\end{proof}

For idempotents with compact difference, the map
$F|_{\operatorname{Ran}E}:\operatorname{Ran}E\to\operatorname{Ran}F$ is Fredholm: $E|_{\operatorname{Ran}F}$ is a two-sided parametrix modulo compact operators. We use
\begin{equation}\label{cl:index}
 \operatorname{ind}(E,F)=\operatorname{ind}\bigl(F|_{\operatorname{Ran}E}\bigr),
 \qquad \nu_{E_0}(E)=\operatorname{ind}(E,E_0).
\end{equation}
The index is locally constant and additive:
\[
 \operatorname{ind}(E,G)=\operatorname{ind}(E,F)+\operatorname{ind}(F,G).
\]
These are the standard Fredholm-pair properties; see, for example, \cite{ref9}. In particular, if $E=F+q$, where $q$ is a finite-rank orthogonal projection annihilated by $F$ on both sides, then $\operatorname{ind}(E,F)=\operatorname{rank}q$. The modification in Lemma~\ref{cl:install} has index zero.

\begin{lemma}[Restricted-unitary path]\label{cl:unitary-path}
If $P,Q$ are orthogonal projections with $P-Q\in L^p$ and $\operatorname{ind}(P,Q)=0$, there is a self-adjoint $K\in L^p$ such that $Q=e^{iK}Pe^{-iK}$.
\end{lemma}

\begin{proof}
Set
\[
 A=QP+(I-Q)(I-P)=I+(Q-P)(2P-I).
\]
Then $AP=QA$, $A-I\in L^p$, and
\[
 A^*A=AA^*=I-(P-Q)^2.
\]
Its kernel is the finite-dimensional space
\[
 M=(PH\cap\ker Q)\oplus(\ker P\cap QH).
\]
The zero-index assumption makes these summands equidimensional. Choose a unitary $W$ on $M$ interchanging them, and let $P_M$ be the orthogonal projection onto $M$. With $B=(A^*A)^{1/2}$, the operators
\[
 V=A(B+P_M)^{-1},\qquad U=V+WP_M
\]
give a unitary $U$ satisfying $UP=QU$. Here $B,P_M$ commute with $P,Q$; $V$ is unitary on $M^\perp$ and zero on $M$. From $(I-B)(I+B)=(P-Q)^2$ we obtain $B-I\in L^p$, and hence $U-I\in L^p$. The spectral theorem for the compact normal operator $U-I$ gives a self-adjoint logarithm $K$ with spectrum in $[-\pi,\pi]$. Since
\[
 |t|\leq (\pi/2)|e^{it}-1|\qquad(-\pi\leq t\leq\pi),
\]
we have $K\in L^p$. Thus $U=e^{iK}$, and $e^{itK}Pe^{-itK}$ is the required $L^p$-continuous path.
\end{proof}

\begin{lemma}[Connected components]\label{cl:components}
For every $n\in\mathbb Z$, the set
\[
 \mathcal C_p^{(n)}(E_0)
 =\{E\in\mathcal C_p(E_0):\nu_{E_0}(E)=n\}
\]
is nonempty and path connected in the $L^p$ topology. These sets are exactly the connected components of $\mathcal C_p(E_0)$.
\end{lemma}

\begin{proof}
Choose a bounded invertible $S$ with $E_0=SPS^{-1}$ for an orthogonal projection $P$. Conjugation by $S$ is a homeomorphism for the $L^p$ topology, although it need not be an $L^p$ isometry. Suppose $E-P\in L^p$, and let $Q=R(E)$ be the orthogonal range projection. The range-projection formula in Lemma~\ref{common:block} and the resolvent identity give $Q-P\in L^p$. Since $QE=E$ and $EQ=Q$, the operators
\[
 Q+t(E-Q),\qquad 0\leq t\leq1,
\]
form an idempotent path of constant range, preserving the relative index. For each $n$, obtain $P_n$ by adding $n$ orthogonal directions to $P$ if $n\geq0$, or removing $|n|$ directions if $n<0$. If $\nu_P(Q)=n$, additivity gives $\operatorname{ind}(P_n,Q)=0$, so Lemma~\ref{cl:unitary-path} joins $P_n$ to $Q$. This proves nonemptiness and path connectedness. Local constancy of the Fredholm index prevents a connected set from meeting distinct labels. Undo the fixed similarity to obtain the claim for $E_0$.
\end{proof}

\begin{lemma}[Fixed zero blocks and square-zero lines]\label{cl:zero-block-lines}
For a finite-rank orthogonal projection $q$ and $n\in\mathbb Z$, let
\[
 \mathcal D_q^{(n)}=
 \{E\in\mathcal C_p(E_0):Eq=qE=0,\ \nu_{E_0}(E)=n\}.
\]
This set is nonempty. Any two of its elements can be joined by a finite chain of complete affine real lines of idempotents, each with square-zero direction and contained in $\mathcal D_q^{(n)}$.
\end{lemma}

\begin{proof}
Apply Lemma~\ref{cl:install} to $V=qH\oplus V_+\oplus V_-$, where $\dim V_\pm\geq|n|$, installing zero on $qH\oplus V_-$ and the identity on $V_+$. The resulting idempotent has index zero. Adding a rank-$n$ subprojection in $V_-$ when $n\geq0$, or removing a rank-$|n|$ subprojection in $V_+$ when $n<0$, produces an element $G_n\in\mathcal D_q^{(n)}$. Every element is $0\oplus E'$ on $qH\oplus q^\perp H$, and $\operatorname{ind}(E',G_n')=0$. Lemma~\ref{cl:components}, applied on $q^\perp H$, therefore gives path connectedness while fixing the zero block.

For sufficiently close idempotents $E,F$, put
\[
 M=E+F-I,\qquad G=EM^{-1}F.
\]
The factorization
\[
 M=(2E-I)\bigl[I+(2E-I)(F-E)\bigr]
\]
gives invertibility when $\|(2E-I)(F-E)\|<1$. The relations $ME=FE=FM$ imply
\[
 EG=G,\quad GE=E,\quad FG=F,\quad GF=G.
\]
Thus $G$ is an idempotent and $(G-E)^2=(F-G)^2=0$. Moreover,
\[
 G-E=EM^{-1}(F-E)(I-E)\in L^p.
\]
If $E,F$ annihilate $q$ on both sides, then $Mq=qM=-q$, so $G$ also annihilates $q$. The square-zero criterion in Lemma~\ref{common:corners} shows that $E+t(G-E)$ and $G+t(F-G)$ are complete lines of idempotents. Their indices are constant by continuity, so both lines remain in $\mathcal D_q^{(n)}$. Along a continuous path $E_s$, the norms $\|2E_s-I\|$ are bounded. A sufficiently fine finite subdivision makes every consecutive pair satisfy the preceding invertibility condition. Replace each short segment by its two lines to obtain the finite chain.
\end{proof}

\begin{theorem}[Minimum-distance equality]\label{cl:minimum}
If $E-F\in L^p$ and $n=\operatorname{ind}(E,F)\ne0$, then
\begin{equation}\label{cl:minimum-bound}
 \|E-F\|_p\geq |n|^{1/p}.
\end{equation}
For $n>0$, equality holds exactly when $E=F+q$, where $q$ is an orthogonal projection of rank $n$ and $Fq=qF=0$. For $n<0$, interchange $E,F$. Consequently,
\begin{equation}\label{cl:component-distance}
 \operatorname{dist}_p\bigl(\mathcal C_p^{(m)}(E_0),\mathcal C_p^{(n)}(E_0)\bigr)
 =|m-n|^{1/p}.
\end{equation}
\end{theorem}

\begin{proof}
For $n>0$, the space $L=\operatorname{Ran}E\cap\ker F$ has dimension at least $n$, and $D=E-F$ is the identity on $L$. The variational principle gives $s_j(D)\geq1$ for $1\leq j\leq n$, proving \eqref{cl:minimum-bound}. Equality forces exactly $n$ nonzero singular values, all equal to $1$. Thus $D$ is a rank-$n$ partial isometry. Its initial space contains $L$, since $Dx=x$ there; the initial space, final space, and $L$ must therefore coincide. Hence $D$ is the orthogonal projection onto $L$. The elementary block observation preceding Lemma~\ref{cl:install} gives $Fq=qF=0$. The converse and the negative case follow immediately.

The lower bound in \eqref{cl:component-distance} follows from index additivity. To attain it, use Lemma~\ref{cl:install} to install equally large finite identity and zero blocks. Adding or removing nested orthogonal directions in this block gives labels $m,n$ with difference a signed orthogonal projection of rank $|m-n|$.
\end{proof}

Equation~\eqref{cl:component-distance} is a distance between sets; it does not assert that every point attains that distance to another component. For example, on $K\oplus K$, the idempotent
\[
 \begin{pmatrix}I&aI\\0&0\end{pmatrix},\qquad a>0,
\]
has no nonzero orthogonal reducing zero or identity block. The equality characterization in Theorem~\ref{cl:minimum} therefore rules out an adjacent minimum pair through this point.

\subsection{Rank-one increments and Wigner normalization}

Let $\mathcal P_1(H)$ denote the rank-one orthogonal projections. Propagation at the $L^1$ endpoint uses the following equality mechanism.

\begin{lemma}[$L^1$ norm rectangle]\label{cl:rectangle}
If $T\in L^1$ and $r,s\in\mathcal P_1(H)$ satisfy
\[
 \|T-r\|_1=\|T\|_1+1,\qquad
 \|T+s\|_1=\|T\|_1+1,\qquad
 \|T+s-r\|_1=\|T\|_1,
\]
then $r=s$ and $Tr=rT=0$.
\end{lemma}

\begin{proof}
Equality in the first triangle inequality gives a contraction $V$ with
\[
 \operatorname{Re}\operatorname{Tr}(V^*T)=\|T\|_1,
 \qquad Vr=V^*r=-r.
\]
Evaluating $V$ on the third equality yields
\[
 \|T\|_1\geq\|T\|_1+1+\operatorname{Re}\operatorname{Tr}(V^*s),
\]
so $Vs=V^*s=-s$. Similarly, a supporting contraction $W$ for $T$ satisfies $Wr=W^*r=r$ and $Ws=W^*s=s$. Write a singular-value expansion
\[
 T=\sum_j\sigma_j u_j\otimes v_j,\qquad \sigma_j>0.
\]
Both contractions attain the $L^1$ norm, so $Vv_j=Wv_j=u_j$. For $z\in rH+sH$,
\[
 0=\langle(V-W)v_j,z\rangle
 =\langle v_j,(V^*-W^*)z\rangle=-2\langle v_j,z\rangle.
\]
The same orthogonality holds for $u_j$. Thus the initial and final supports of $T$ lie in $(rH+sH)^\perp$, whereas $s-r$ is supported on $rH+sH$. Their $L^1$ norms add, and the third equality forces $s=r$. Support orthogonality then gives $Tr=rT=0$.
\end{proof}

\begin{proposition}[Global orthogonal increments]\label{cl:increments}
For a surjective isometry $\Phi$ as in Theorem~\ref{cl:main}, there are $\varepsilon\in\{1,-1\}$ and a bijection $r:\mathcal P_1(H)\to\mathcal P_1(H)$ such that
\begin{equation}\label{cl:increment-formula}
 \Phi(E+q)-\Phi(E)=\varepsilon r(q)\qquad(Eq=qE=0),
\end{equation}
and $\operatorname{Tr}(r(q)r(s))=\operatorname{Tr}(qs)$ for all $q,s$.
\end{proposition}

\begin{proof}
The isometry permutes the connected components and, by \eqref{cl:component-distance}, acts on their integer labels as $n\mapsto\varepsilon n+c$. The minimum-distance equality gives
\[
 \Phi(E+q)-\Phi(E)=\varepsilon r_E(q)
\]
with $r_E(q)$ rank-one orthogonal.

Fix $q$ and let $E_t=E+tN$ be a complete line from Lemma~\ref{cl:zero-block-lines}. If $p>1$, strict convexity of $L^p$ makes both real isometric lines $\Phi(E_t)$ and $\Phi(E_t+q)$ affine. Their distance is constantly $1$, so their slopes agree and the increment is independent of $t$. For $p=1$, the four points $E_s,E_t,E_s+q,E_t+q$ have cross-distances $|t-s|\|N\|_1+1$ and corresponding distances $|t-s|\|N\|_1$; the zero block makes the relevant supports orthogonal. Apply Lemma~\ref{cl:rectangle} to their image differences, including the orientation $\varepsilon$, to get $r_{E_s}(q)=r_{E_t}(q)$. The finite-chain lemma therefore makes this projection depend only on the index label; write it as $r_n(q)$.

For orthogonal $q,s$, install a common zero block at any label $n$. Computing $\Phi(E+q+s)-\Phi(E)$ in the two orders gives
\[
 a+c=b+d,\qquad
 a=r_n(q),\quad b=r_n(s),\quad c=r_{n+1}(s),\quad d=r_{n+1}(q).
\]
Since $\|q-s\|_p=2^{1/p}$, the same holds for $a,b$, hence $ab=0$. Multiply the identity by $a$ and take traces:
\[
 1+\operatorname{Tr}(ac)=\operatorname{Tr}(ad)\leq1.
\]
Thus $d=a$, removing the dependence on $n$.

For arbitrary $q,s$, choose a common zero block and compare $E+q,E+s$. The identity
\[
 \|q-s\|_p^p=2\bigl(1-\operatorname{Tr}(qs)\bigr)^{p/2}
\]
gives preservation of transition probabilities and injectivity of $r$. To prove surjectivity, fix a target rank-one projection $h$. If $\varepsilon=1$, install a target zero block $Yh=hY=0$ and compare $Y,Y+h$. Their preimages have index difference $1$ and distance $1$, hence are $X,X+q$ by Theorem~\ref{cl:minimum}; thus $r(q)=h$. If $\varepsilon=-1$, install a target identity block $Yh=hY=h$ and compare $Y,Y-h$. Reversal of the index orientation makes their ordered preimages again $X,X+q$, with image difference $-h$. Equation~\eqref{cl:increment-formula} again gives $r(q)=h$.
\end{proof}

Wigner's theorem \cite{ref17} supplies a unitary or antiunitary $U$ such that $r(q)=UqU^*$. Conjugating the output by $U^*$ and, if necessary, complementing it gives a surjective $L^p$ isometry $\Psi$ between possibly different perturbation classes such that
\begin{equation}\label{cl:normalized-increments}
 \Psi(E+q)=\Psi(E)+q\qquad(Eq=qE=0).
\end{equation}
Fix orthogonal unit vectors $x,y$, and put
\[
 N=x\otimes y,\qquad q_x=x\otimes x,\qquad q_y=y\otimes y,\qquad P=q_x+q_y.
\]
Use Lemma~\ref{cl:install} to choose $G$ reduced by $PH$, with $G|_{PH}=q_x$. Applying \eqref{cl:normalized-increments} to $G-q_x,G,G+q_y$ shows that $PH$ also reduces $F=\Psi(G)$ and $F|_{PH}=q_x$. Define
\[
 f(z)=\Psi(G+zN)-\Psi(G).
\]
Then
\begin{equation}\label{cl:plane-probes}
 \begin{aligned}
 f(0)&=0,& \|f(z)-f(w)\|_p&=|z-w|,\\
 \|f(z)-(h-q_x)\|_p&=\|zN-(h-q_x)\|_p
 &&\bigl(h\in\mathcal P_1(H),\ h\leq P\bigr).
 \end{aligned}
\end{equation}
The second equality follows from $\Psi(G-q_x+h)=\Psi(G)-q_x+h$. These projection probes supply the local information required below.

\subsection{\texorpdfstring{Local plane rigidity for $1<p<\infty$}{Local plane rigidity for 1 < p < infinity}}

\begin{lemma}[Strictly convex local rigidity]\label{cl:local-strict}
For $1<p<\infty$, the function in \eqref{cl:plane-probes} has one of the forms
\begin{equation}\label{cl:local-alternatives}
 f(z)=zN\quad(z\in\mathbb C),\qquad
 f(z)=\overline zN^*\quad(z\in\mathbb C).
\end{equation}
\end{lemma}

\begin{proof}
Strict convexity gives uniqueness of midpoints in an isometric image of $\mathbb C$, so $f$ is real linear. For $B=f(z)$ the whole line $F+tB$ consists of idempotents. Comparison of coefficients gives
\begin{equation}\label{cl:local-square-zero}
 FB+BF=B,\qquad B^2=0,\qquad \|B\|_p=|z|.
\end{equation}
For rank-one $h\leq P$ with $h\ne q_x$, put $D=h-q_x$. Then $D$ is self-adjoint and $D^2=\rho^2P$ for some $\rho>0$, and
\[
 D\|\cdot\|_p(D)[K]
 =\frac{\operatorname{Re}\operatorname{Tr}(DK)}{2^{(p-1)/p}\rho}.
\]
Differentiate the probe identity at $tz$, $t=0$. Since the matrices $h-q_x$ span the real space of traceless self-adjoint matrices on $PH$, we obtain
\[
 P(B+B^*)P=zN+\overline zN^*+2aP
\]
for a real $a$. The first relation in \eqref{cl:local-square-zero}, compressed to $PH$ where $F=q_x$, makes both diagonal coefficients of $PBP$ zero, so $a=0$.

The initial and final support spaces of a square-zero operator are orthogonal. The nonzero singular values of $B+B^*$ are therefore those of $B$, each repeated twice. Consequently,
\[
 \|B+B^*\|_p^p=2\|B\|_p^p=2|z|^p
 =\|P(B+B^*)P\|_p^p.
\]
To analyze this compression equality, set $T=B+B^*$ and $R=2P-I$. Pinching gives
\[
 \|PTP\|_p\leq\left\|\frac{T+RTR}{2}\right\|_p
 \leq\|T\|_p=\|PTP\|_p.
\]
Equality in the second inequality and strict convexity force $T=RTR$. The direct-sum formula in the first equality then gives $P^\perp TP^\perp=0$. Hence
\begin{equation}\label{cl:local-symmetric-part}
 B+B^*=zN+\overline zN^*.
\end{equation}
Write $B=X+iY$ with $X,Y$ self-adjoint. The equation $B^2=0$ gives $X^2=Y^2$ and $XY+YX=0$, while \eqref{cl:local-symmetric-part} gives $X^2=|z|^2P/4$. Thus $Y$ annihilates $P^\perp H$ and $B=PBP$. The first relation in \eqref{cl:local-square-zero} now gives
\[
 B=\begin{pmatrix}0&b\\c&0\end{pmatrix},\qquad bc=0,\qquad b+\overline c=z.
\]
Therefore $B=zN$ or $B=\overline zN^*$. For $z\ne0$ the alternatives are separated by $2^{1/p}|z|$; continuity and connectedness of the punctured plane make the choice uniform.
\end{proof}

\subsection{\texorpdfstring{Local rigidity at the $L^1$ endpoint}{Local rigidity at the L1 endpoint}}

The endpoint argument cannot assume that images of affine lines remain affine. We give the required analytic ingredients explicitly.

\begin{lemma}[Quantitative $L^1$ compression]\label{cl:compression}
If $T\in L^1$, $\|T\|_1=1$, and $P$ is an orthogonal projection, then, with $m=\|PTP\|_1$,
\[
 \|T-PTP\|_1\leq2\sqrt{1-m^2}.
\]
\end{lemma}

\begin{proof}
For the polar decomposition $T=V|T|$, set $A=V|T|^{1/2}$ and $B=|T|^{1/2}$. Then $T=AB$ and $\|A\|_2=\|B\|_2=1$. With $\alpha=\|PA\|_2^2$ and $\beta=\|BP\|_2^2$, the $L^2\cdot L^2\subset L^1$ H\"older inequality gives $m\leq\sqrt{\alpha\beta}$. Both numbers are at most $1$, so each is at least $m^2$. The decomposition
\[
 T-PTP=P^\perp AB+PABP^\perp
\]
therefore yields
\[
 \|T-PTP\|_1\leq\sqrt{1-\alpha}+\sqrt{1-\beta}
 \leq2\sqrt{1-m^2}.
\]
\end{proof}

\begin{lemma}[Two-dimensional $L^1$ norm rigidity]\label{cl:two-dimensional}
On $\mathbb C^2$, let
\[
 N=\begin{pmatrix}0&1\\0&0\end{pmatrix},\qquad
 H_c=cN+\overline cN^*.
\]
If
\[
 \|Y-H_c\|_1=|1-c|+|c|\qquad(c\in\mathbb C),
\]
then $Y=N$ or $Y=N^*$.
\end{lemma}

\begin{proof}
Put $J=N+N^*$. The choices $c=0,1$ give $\|Y\|_1=\|J-Y\|_1=1$, whereas $\|J\|_1=2$. Since $J$ is a self-adjoint unitary, equality in trace duality makes $B=JY$ and $I-B=J(J-Y)$ positive. Thus
\[
 B=\begin{pmatrix}u&v\\\overline v&1-u\end{pmatrix},\qquad
 0\leq u\leq1,\quad |v|^2\leq u(1-u),\qquad
 Y=\begin{pmatrix}\overline v&1-u\\u&v\end{pmatrix}.
\]
Set $d=u(1-u)-|v|^2\geq0$. The identity
\[
 \|M\|_1^2=\operatorname{Tr}(M^*M)+2|\det M|
\]
for $2\times2$ matrices, with $c=it$, $t>0$, gives
\[
 \|Y-H_{it}\|_1^2
 =1-2d+2t^2+2\sqrt{(d+t^2)^2+t^2(1-2u)^2}.
\]
The hypothesis consequently implies
\[
 \sqrt{(d+t^2)^2+t^2(1-2u)^2}=d+t\sqrt{1+t^2}.
\]
If $d>0$, the right derivatives at zero are respectively $0$ and $1$, a contradiction. If $d=0$, division by $t$ and squaring give $(1-2u)^2=1$. Hence $u\in\{0,1\}$ and $v=0$, proving the assertion.
\end{proof}

\begin{lemma}[One-sided $L^1$ norm derivative]\label{cl:derivative}
Let $A=W|A|\in L^1$ be a polar decomposition, and put $P_A=W^*W$, $Q_A=WW^*$. For $K\in L^1$,
\begin{equation}\label{cl:derivative-formula}
 D^+\|\cdot\|_1(A)[K]
 =\operatorname{Re}\operatorname{Tr}(W^*K)
 +\|(I-Q_A)K(I-P_A)\|_1.
\end{equation}
\end{lemma}

\begin{proof}
If a contraction $X$ supports the norm at $A$, expansion in the singular vectors of $A$ shows that $Xv_j=Wv_j$ and $X^*Wv_j=v_j$ on every nonzero singular vector. All its off-diagonal support blocks therefore vanish, and
\[
 X=W+Z,\qquad Z=(I-Q_A)Z(I-P_A),\qquad\|Z\|\leq1.
\]
Conversely, every such $X$ is a supporting contraction, because the two summands have orthogonal initial and final supports. The right directional derivative of a continuous convex norm is the supremum of its supporting real-linear functionals. Taking this supremum and using $L^1$ duality on the complementary supports gives \eqref{cl:derivative-formula}.
\end{proof}

\begin{lemma}[Endpoint local rigidity]\label{cl:local-endpoint}
For $p=1$, the function in \eqref{cl:plane-probes} also has one of the forms \eqref{cl:local-alternatives}.
\end{lemma}

\begin{proof}
For $z\ne0$, put $t=|z|$, $\lambda=z/|z|$, and $X_z=f(z)/t$, so $\|X_z\|_1=1$. Let $h_{t,c}$ be the rank-one orthogonal projection onto $\mathbb C(x+t\overline c y)$. Its exact graph-projection matrix gives
\[
 \frac{h_{t,c}-q_x}{t}=H_c+O(t),\qquad H_c=cN+\overline cN^*.
\]
The error is in $L^1$ norm and is uniform for bounded $c$. The probe identity yields
\begin{equation}\label{cl:endpoint-scaled-probe}
 \|X_z-H_c\|_1=\|\lambda N-H_c\|_1+O(t).
\end{equation}
For $c=2\lambda$, put $J_\lambda=\lambda N+\overline\lambda N^*$. It is a self-adjoint unitary on $PH$, zero outside it, and $\|\lambda N-2J_\lambda\|_1=3$. Duality with $-J_\lambda$ gives
\[
 3+O(t)=\|X_z-2J_\lambda\|_1
 \geq4-\operatorname{Re}\operatorname{Tr}(J_\lambda X_z).
\]
Thus $\|PX_zP\|_1\geq1-Ct$, with $C$ independent of the phase of $z$. Lemma~\ref{cl:compression} implies
\begin{equation}\label{cl:endpoint-compactness}
 \|X_z-PX_zP\|_1=O(|z|^{1/2})\longrightarrow0.
\end{equation}
This estimate supplies the required compactness; the unit ball of $L^1$ is not norm compact.

For a sequence $z_j\to0$, $z_j\ne0$, pass to a subsequence with $\lambda_j\to\lambda$. The finite-dimensional compressions and \eqref{cl:endpoint-compactness} permit a further subsequence for which $X_{z_j}\to X\in P\mathcal B(H)P$ in $L^1$ norm. Passing to the limit in \eqref{cl:endpoint-scaled-probe} gives
\[
 \|X-H_c\|_1=|\lambda-c|+|c|\qquad(c\in\mathbb C).
\]
Conjugate by $D_\lambda=\operatorname{diag}(\lambda,1)$ on $PH$. Since
\[
 D_\lambda^*\lambda ND_\lambda=N,\qquad
 D_\lambda^*H_cD_\lambda=H_{\overline\lambda c},
\]
Lemma~\ref{cl:two-dimensional} gives $X=\lambda N$ or $X=\overline\lambda N^*$. Hence
\[
 \operatorname{dist}_1\left(X_z,
 \left\{\frac{z}{|z|}N,\frac{\overline z}{|z|}N^*\right\}\right)
 \longrightarrow0.
\]
The candidates are always at distance $2$. On a sufficiently small punctured disk, continuity gives two disjoint open alternatives; connectedness makes one alternative uniform. We have proved either
\begin{equation}\label{cl:endpoint-first-order}
 f(z)=zN+o(|z|),\qquad\text{or}\qquad
 f(z)=\overline zN^*+o(|z|).
\end{equation}
In particular, real Fr\'echet differentiability at the origin follows from the metric probes; it has not been assumed.

Assume the first case. Fix $w=r\lambda\ne0$, where $r=|w|$, and write $A=f(w)=W|A|$, with support projections $P_A,Q_A$. The isometry and \eqref{cl:endpoint-first-order} imply, for every $\zeta\in\mathbb C$,
\[
 D^+\|\cdot\|_1(A)[-\zeta N]=-\operatorname{Re}(\overline\lambda\zeta).
\]
Writing $\beta=\|(I-Q_A)N(I-P_A)\|_1$, the derivative formula becomes
\[
 -\operatorname{Re}\{\zeta\operatorname{Tr}(W^*N)\}+|\zeta|\beta
 =-\operatorname{Re}(\overline\lambda\zeta).
\]
The equations for $\zeta$ and $-\zeta$ give $\beta=0$, and variation of the phase gives $\operatorname{Tr}(W^*N)=\overline\lambda$. Equality for this rank-one coefficient of a contraction yields
\begin{equation}\label{cl:endpoint-polar-values}
 Wy=\lambda x,\qquad W^*x=\overline\lambda y,\qquad
 P_Ay=y,\quad Q_Ax=x.
\end{equation}
Next differentiate the projection-probe equality with $h_{t,c}$ at $t=0+$. It gives
\[
 D^+\|\cdot\|_1(A)[-H_c]
 =D^+\|\cdot\|_1(wN)[-H_c].
\]
Using \eqref{cl:endpoint-polar-values} and putting
\[
 \eta=\operatorname{Tr}(W^*N^*),\qquad
 \alpha=\|(I-Q_A)y\|\,\|(I-P_A)x\|,
\]
this reduces to
\[
 -\operatorname{Re}(\overline c\eta)+|c|\alpha=|c|.
\]
Add the equations for $c$ and $-c$ to obtain $\alpha=1$. Each factor is at most $1$, so
\begin{equation}\label{cl:endpoint-support-values}
 P_Ay=y,\quad P_Ax=0,\quad Q_Ax=x,\quad Q_Ay=0.
\end{equation}

Set $a=\langle|A|y,y\rangle$, so $0\leq a\leq r=\operatorname{Tr}|A|$. We claim that
\begin{equation}\label{cl:endpoint-block-norms}
 PAP=a\lambda N,\qquad \|P^\perp AP^\perp\|_1=r-a.
\end{equation}
Indeed, $P_Ax=0$ gives $|A|x=0$. The positive compression
\[
 K=(I-q_y)|A|(I-q_y)
\]
satisfies $K=P_AKP_A$ and $\operatorname{Tr}K=r-a$. Relations \eqref{cl:endpoint-polar-values}--\eqref{cl:endpoint-support-values} show that $W$ maps $P_AH\cap y^\perp$ isometrically into $P^\perp H$. Relative to $\mathbb Cx\oplus\mathbb Cy\oplus P^\perp H$, write
\[
 |A|=\begin{pmatrix}0&0&0\\0&a&b^*\\0&b&D\end{pmatrix},
 \qquad K=0\oplus0\oplus D\geq0.
\]
Direct compression of $W|A|$ gives $P^\perp AP^\perp=WK$. Since $W^*WK=K$, its $L^1$ norm is $\operatorname{Tr}K=r-a$. This proves \eqref{cl:endpoint-block-norms}; no vanishing of the mixed blocks of $A$ has been asserted. Pinching therefore gives, for every $D_0$ supported on $PH$,
\begin{equation}\label{cl:endpoint-pinching}
 \|A-D_0\|_1\geq\|a\lambda N-D_0\|_1+r-a.
\end{equation}

For $s\geq0$, let $h_s$ project onto $\mathbb C(x+s\overline\lambda y)$ and put $D_s=h_s-q_x$. A direct $2\times2$ calculation gives, for $u\geq0$,
\begin{equation}\label{cl:endpoint-graph-norm}
 \|u\lambda N-D_s\|_1=
 \begin{cases}
 u,&0\leq s\leq u,\\[2pt]
 \displaystyle\sqrt{u^2+\frac{4s(s-u)}{1+s^2}},&s>u.
 \end{cases}
\end{equation}
For completeness, the graph-projection matrix yields
\[
 \operatorname{Tr}\bigl((u\lambda N-D_s)^*(u\lambda N-D_s)\bigr)
 =u^2+\frac{2s(s-u)}{1+s^2},\qquad
 \det(u\lambda N-D_s)=\frac{s(u-s)}{1+s^2};
\]
the $L^1$ determinant identity proves \eqref{cl:endpoint-graph-norm}. If $a<r$, choose $a<s<r$. The probe equality gives
\[
 \|A-D_s\|_1=\|r\lambda N-D_s\|_1=r,
\]
whereas \eqref{cl:endpoint-pinching}--\eqref{cl:endpoint-graph-norm} give a value strictly greater than $r$. Hence $a=r$. Positivity forces $|A|=r(y\otimes y)$, and \eqref{cl:endpoint-polar-values} gives $A=wN$. This holds for every $w$, proving $f(z)=zN$. In the second case of \eqref{cl:endpoint-first-order}, apply the same argument to $f(z)^*$, whose probe norms agree because $h-q_x$ is self-adjoint, to obtain $f(z)=\overline zN^*$.
\end{proof}

\subsection{Coherence and recovery of the entire class}

\begin{lemma}[Uniform local alternative]\label{cl:coherence}
The choice in \eqref{cl:local-alternatives} is independent of the compatible base point $G$ and of the unit rank-one square-zero direction $N$.
\end{lemma}

\begin{proof}
Fix $N$. If compatible $G_1,G_2$ had different choices, the difference of their image points at $G_j+tN$, for real $t$, would be
\[
 \Psi(G_1)-\Psi(G_2)+t(N-N^*).
\]
Its $L^p$ norm must remain $\|G_1-G_2\|_p$, contradicting the reverse triangle inequality as $|t|\to\infty$.

To compare directions, the space
\[
 \mathcal N_1=\{x\otimes y:x\perp y,\ \|x\|=\|y\|=1\}
\]
is path connected by transitivity of the unitary action and the existence of unitary paths. We verify the local continuity needed for the $L^p$ topology. Fix $N_0=x_0\otimes y_0$ and $P_0=q_{x_0}+q_{y_0}$. For $N$ near $N_0$, put
\[
 x_N=\frac{NN^*x_0}{\|NN^*x_0\|},\qquad
 y_N=N^*x_N,\qquad P_N=q_{x_N}+q_{y_N},
\]
and
\[
 Z_N=x_N\otimes x_0+y_N\otimes y_0+(I-P_N)(I-P_0).
\]
These operators are norm continuous, $Z_{N_0}=I$, and
\[
 Z_N^*Z_N=P_0+(I-P_0)(I-P_N)(I-P_0).
\]
Near $N_0$, the polar unitary $V_N=Z_N(Z_N^*Z_N)^{-1/2}$ is norm continuous and sends $x_0,y_0$ to $x_N,y_N$; hence $V_NN_0V_N^*=N$. Moreover,
\[
 \operatorname{rank}(Z_N-I)\leq4,\qquad
 \operatorname{rank}(Z_N^*Z_N-I)\leq8,\qquad
 \operatorname{rank}(V_N-I)\leq12,
\]
the last bound following by functional calculus. Choose a compatible $G_0$, and put $G_N=V_NG_0V_N^*$. These operators belong to the same perturbation class and are compatible with $N$. The uniformly bounded ranks of $V_N-V_{N'}$ give
\[
 \|G_N-G_{N'}\|_p
 \leq2\,24^{1/p}\|G_0\|\,\|V_N-V_{N'}\|.
\]
Thus $N\mapsto\Psi(G_N+N)-\Psi(G_N)$ is $L^p$ continuous. Its value is either $N$ or $N^*$, separated by distance $2^{1/p}$. The choice is locally constant, and path connectedness makes it global.
\end{proof}

After optional adjunction of $\Psi$, which preserves \eqref{cl:normalized-increments}, we may therefore arrange
\begin{equation}\label{cl:normalized-lines}
 \Psi(G+zN)=\Psi(G)+zN
\end{equation}
for every compatible $G$, every $N\in\mathcal N_1$, and $z\in\mathbb C$.

\begin{lemma}[Rank-one probes at infinity]\label{cl:infinity-probes}
For $N=x\otimes y\in\mathcal N_1$ and $A\in L^p$,
\begin{align}
 \|te^{i\theta}N-A\|_p
 &=t-\operatorname{Re}\{e^{-i\theta}\operatorname{Tr}(N^*A)\}+o(1),
 &&1<p<\infty,\label{cl:infinity-strict}\\
 \|te^{i\theta}N-A\|_1
 &=t-\operatorname{Re}\{e^{-i\theta}\operatorname{Tr}(N^*A)\}
 +\|(I-q_x)A(I-q_y)\|_1+o(1).
 \label{cl:infinity-endpoint}
\end{align}
If $K\in L^p$ satisfies $\operatorname{Tr}(N^*K)=0$ for every $N\in\mathcal N_1$, then $K=0$.
\end{lemma}

\begin{proof}
Divide each norm by $t$. For $p>1$, Fr\'echet differentiability at $e^{i\theta}N$ gives \eqref{cl:infinity-strict}; for $p=1$, use Lemma~\ref{cl:derivative}. The complementary support projections are $I-q_y$ initially and $I-q_x$ finally, giving \eqref{cl:infinity-endpoint}. For the last assertion, the coefficients say that $Ky\in\mathbb Cy$ for every nonzero $y$. Apply this to independent vectors and their sum to obtain $K=\lambda I$. Compactness and infinite dimensionality force $\lambda=0$.
\end{proof}

\begin{proof}[Proof of Theorem~\ref{cl:main}]
Proposition~\ref{cl:increments}, Wigner normalization, the two local rigidity lemmas, and Lemma~\ref{cl:coherence} give \eqref{cl:normalized-lines}. For any $E$ in the source class and $N\in\mathcal N_1$, choose a compatible $G$ using Lemma~\ref{cl:install}. Isometry gives
\[
 \|te^{i\theta}N-(E-G)\|_p
 =\|te^{i\theta}N-(\Psi(E)-\Psi(G))\|_p.
\]
For $p>1$, \eqref{cl:infinity-strict} recovers the real part of the $N$ coefficient. For $p=1$, subtract the equalities for $\theta$ and $\theta+\pi$ before taking the limit; the phase-independent complementary $L^1$ terms cancel. Varying $\theta$ gives
\[
 \operatorname{Tr}(N^*(E-G))
 =\operatorname{Tr}(N^*(\Psi(E)-\Psi(G))).
\]
Use the same $G,N$ for $F$ and subtract. The final assertion of Lemma~\ref{cl:infinity-probes} yields
\[
 \Psi(E)-\Psi(F)=E-F.
\]
Hence $\Psi(E)=E+C$ for one fixed bounded operator $C$. This does not assert $C\in L^p$, since the source and normalized target classes may differ.

Both $G+zN$ and $G+C+zN$ are idempotents. Subtracting the coefficients of $z$ in their idempotent equations gives $CN+NC=0$ for every $N=x\otimes y\in\mathcal N_1$. Fix $y$ and a vector $z\perp y$; evaluating this identity at $z$ gives $\langle Cz,y\rangle x=0$. Thus $C^*y\in\mathbb Cy$ for every $y$. Linearity on independent vectors makes $C^*$ scalar, and $CN+NC=0$ forces that scalar to vanish. The normalized map is the identity. Undo adjunction, Wigner conjugation, and complementation to obtain \eqref{cl:standard-forms}. The converse follows from invariance of the $L^p$ norm under these operations and the stated base-point compatibility condition.
\end{proof}

\subsection{Finite-type sectors and classification of the full space}

Set
\[
 \begin{aligned}
 \mathcal F&=\{E\in\mathcal I(H):\operatorname{rank}E<\infty\},\\
 \mathcal F^\vee&=\{E\in\mathcal I(H):\operatorname{rank}(I-E)<\infty\},\\
 \mathcal J&=\{E\in\mathcal I(H):\operatorname{rank}E=\operatorname{rank}(I-E)=\infty\}.
 \end{aligned}
\]
These three sets are disjoint and exhaust $\mathcal I(H)$.

\begin{lemma}[Complete finite-rank probes]\label{cl:finite-probes}
If a surjective $d_p$ isometry $\Theta:\mathcal F\to\mathcal F$ fixes every rank-one idempotent, it fixes every finite-rank idempotent.
\end{lemma}

\begin{proof}
For $E\in\mathcal F$, put $A=\Theta(E)$. Given orthogonal unit vectors $x,y$, let $q=x\otimes x$ and $N=x\otimes y$. Every $R_{t,\theta}=q+te^{i\theta}N$ is a rank-one idempotent, hence fixed. Therefore
\[
 \|te^{i\theta}N-(A-q)\|_p
 =\|te^{i\theta}N-(E-q)\|_p.
\]
Lemma~\ref{cl:infinity-probes}, using opposite phases for $p=1$, gives $\operatorname{Tr}(N^*(A-E))=0$ for all such pairs. Its separation assertion gives $A=E$.
\end{proof}

\begin{proposition}[One implementation for all finite ranks]\label{cl:finite-sector}
Every surjective $d_p$ isometry of $\mathcal F$ has the form $E\mapsto UEU^*$ or $E\mapsto UE^*U^*$, with one fixed unitary or antiunitary $U$ and one fixed choice of adjunction.
\end{proposition}

\begin{proof}
The rank strata $\mathcal I_k(H)$, $k\geq0$, are the connected components of $\mathcal F$. Indeed, $R(E)+t(E-R(E))$ joins a finite-rank idempotent to its range projection, and fixed-rank orthogonal projections are joined by finite-rank unitary conjugacies. Rank is locally constant on idempotents. The dimension argument of Theorem~\ref{cl:minimum} and nested orthogonal projections give
\[
 \operatorname{dist}_p\bigl(\mathcal I_k(H),\mathcal I_l(H)\bigr)=|k-l|^{1/p}.
\]
The point $0$ is the unique isolated point of $\mathcal F$, so it is fixed. Distances to its component identify every rank stratum, and $\Theta$ preserves $\mathcal I_1(H)$. Apply the fixed-rank $L^p$ classification from Section~\ref{sec:fixed-rank} at $k=1$, which is valid for every $1\leq p<\infty$, and normalize by its standard implementation. The resulting map fixes every rank-one idempotent, so Lemma~\ref{cl:finite-probes} makes it the identity on $\mathcal F$. Thus one normalization from the fixed-rank theorem, together with the complete probes, enforces a common implementation over all finite ranks.
\end{proof}

\begin{theorem}[The full extended $L^p$ metric]\label{cl:full-space}
For $1\leq p<\infty$, a surjective isometry of $(\mathcal I(H),d_p)$ preserves or interchanges $\mathcal F$ and $\mathcal F^\vee$, and preserves $\mathcal J$. Its restriction to $\mathcal F$ has the implementation in Proposition~\ref{cl:finite-sector}, followed by complementation if the finite-type sectors are interchanged. After identifying $\mathcal F^\vee$ with $\mathcal F$ by complementation, its restriction there has a second, independent implementation of the same kind. On each finite-distance class in $\mathcal J$, the restriction has a form in \eqref{cl:standard-forms}; the images give a bijective matching of classes. Conversely, every compatible piecewise construction of this kind is a surjective $d_p$ isometry.
\end{theorem}

\begin{proof}
The only isolated points of the extended metric space are $0$ and $I$. For $E\ne0,I$, choose unit vectors $x\in\operatorname{Ran}E$ and $y\in(\operatorname{Ran}E)^\perp$. Then $N=x\otimes y$ satisfies $EN=N$, $NE=0$, $N^2=0$, so $E+tN$ are distinct idempotents with $d_p(E,E+tN)=|t|$. On the other hand, every nonzero idempotent has operator norm at least $1$; hence $0$ is isolated, and complementation gives the same conclusion for $I$.

The finite-distance class of $0$ is $\mathcal F$, because a compact idempotent has finite rank; the class of $I$ is $\mathcal F^\vee$. Preservation of the unordered pair $\{0,I\}$ yields the asserted partition. After optional global complementation, apply Proposition~\ref{cl:finite-sector} on $\mathcal F$ and, using complementation, independently on $\mathcal F^\vee$. Every finite-distance class in $\mathcal J$ is some $\mathcal C_p(E_0)$, and Theorem~\ref{cl:main} classifies its restricted isometry.

Conversely, the prescribed maps are bijective isometries on each finite-distance class. Distinct classes have distance $+\infty$ between every pair of points, both before and after a bijective matching. Their independent assembly therefore preserves the extended metric.
\end{proof}

\subsection{Uniqueness, group structure, and class invariants}

Let $c(E)=I-E$. Call the maps $U\tau(\cdot)U^*$, with $\tau\in\{\mathrm{id},*,c,c*\}$, \emph{standard transformations}.

\begin{proposition}[Uniqueness on a complete class]\label{cl:uniqueness}
Suppose that, on an infinite-rank, infinite-corank $L^p$ perturbation class,
\[
 U\tau(E)U^*=V\rho(E)V^*,
\]
where $U,V$ are unitary or antiunitary and $\tau,\rho\in\{\mathrm{id},*,c,c*\}$. Then $\tau=\rho$ and $U=\lambda V$ for some $\lambda\in\mathbb T$.
\end{proposition}

\begin{proof}
Given arbitrary $q\in\mathcal P_1(H)$, install $E$ with $Eq=qE=0$ and compare the formulas at $E,E+q$. The complement signs agree, since a nonzero orthogonal projection cannot equal the negative of another such projection, and $UqU^*=VqV^*$ for every $q$. Thus $W=V^*U$ fixes every complex line. In an orthonormal basis, write $We_j=\lambda_j e_j$. Invariance of $\mathbb C(e_i+e_j)$ makes every $\lambda_j$ equal. A conjugate-linear $W$ would then fail to preserve $\mathbb C(e_i+ie_j)$, so $W$ is a scalar linear operator. This proves projective uniqueness and equality of the implementers' linearity types. After removing the common conjugation and complement sign, a remaining difference between $\tau,\rho$ would be adjunction, forcing every element of the class to be self-adjoint. Installing a two-dimensional block $q_x$ and adding $x\otimes y$ gives a contradiction.
\end{proof}

In particular, if $\operatorname{PU}^{\pm}(H)$ denotes the unitary--antiunitary group modulo scalar phases, then
\[
 G_f:=\operatorname{Iso}(\mathcal F,d_p)
 \cong\operatorname{PU}^{\pm}(H)\times\mathbb Z_2,
 \qquad
 G\cong\operatorname{PU}^{\pm}(H)\times\mathbb Z_2\times\mathbb Z_2
\]
is the full standard transformation group, with the extra factors representing adjunction and complementation. These operations commute on operators. Writing $G_\infty=\operatorname{Iso}(\mathcal J,d_p)$, the extended-metric partition gives
\begin{equation}\label{cl:full-group}
 \operatorname{Iso}(\mathcal I(H),d_p)
 \cong\bigl((G_f\times G_f)\rtimes\mathbb Z_2\bigr)\times G_\infty,
\end{equation}
where the semidirect factor interchanges the two finite-type implementations.

\begin{proposition}[Exact matching criterion]\label{cl:matching}
Two infinite-rank, infinite-corank $L^p$ perturbation classes admit a surjective metric isometry exactly when some unitary or antiunitary $U$ and $\tau\in\{\mathrm{id},*,c,c*\}$ satisfy
\begin{equation}\label{cl:matching-condition}
 F_0-U\tau(E_0)U^*\in L^p.
\end{equation}
\end{proposition}

\begin{proof}
Every standard transformation $T$ sends $\mathcal C_p(E_0)$ onto $\mathcal C_p(T(E_0))$, so \eqref{cl:matching-condition} is sufficient. Necessity follows by evaluating Theorem~\ref{cl:main} at $E_0$.
\end{proof}

For completeness, the group on the infinite sector has an unrestricted wreath-product description. Let $\mathfrak C_p$ be the set of infinite-rank, infinite-corank finite-distance classes and $\mathfrak O_p=\mathfrak C_p/G$. For each orbit $\mathcal O$, choose a representative class $C_{\mathcal O}$ and put
\[
 H_{\mathcal O}=\{g\in G:g(C_{\mathcal O})=C_{\mathcal O}\}.
\]
For every $C\in\mathcal O$, choose a standard transformation $g_C$ with $g_C(C_{\mathcal O})=C$. Theorem~\ref{cl:main} and uniqueness identify each isometry of the union of these classes with a permutation $\sigma\in\operatorname{Sym}(\mathcal O)$ and an arbitrary family
\[
 (h_C)_{C\in\mathcal O}\in H_{\mathcal O}^{\mathcal O},
 \qquad \Phi|_C=g_{\sigma(C)}\circ h_C\circ g_C^{-1}.
\]
Consequently,
\begin{equation}\label{cl:wreath-product}
 G_\infty\cong
 \prod_{\mathcal O\in\mathfrak O_p}
 \bigl(H_{\mathcal O}^{\mathcal O}\rtimes\operatorname{Sym}(\mathcal O)\bigr).
\end{equation}
There is no finite-support restriction. This description depends on the chosen representatives and identifications.

\begin{proposition}[An essential-spectral invariant]\label{cl:essential-spectrum}
The set
\[
 \Sigma(E)=\sigma_{\mathrm{ess}}(|E-E^*|)
\]
is constant on each $L^p$ perturbation class and is preserved by every standard transformation.
\end{proposition}

\begin{proof}
If $E-F\in L^p$, the skew-adjoint operators $E-E^*$ and $F-F^*$ have the same Calkin image. Continuous functional calculus in the Calkin algebra gives identical images of their absolute values, hence the same essential spectrum. Adjunction and complementation change the skew-adjoint difference only by a sign. Unitary or antiunitary conjugation preserves the essential spectrum of its positive absolute value. These facts prove both assertions.
\end{proof}

\begin{corollary}[Continuum many metric types]\label{cl:continuum}
For every $1\leq p<\infty$, the infinite-rank, infinite-corank $L^p$ perturbation classes have exactly $2^{\aleph_0}$ standard orbits and exactly $2^{\aleph_0}$ surjective metric-isometry types.
\end{corollary}

\begin{proof}
On $H=\ell^2\oplus\ell^2$, take
\[
 E_a=\begin{pmatrix}I&aI\\0&0\end{pmatrix}\qquad(a\geq0).
\]
Then $|E_a-E_a^*|=aI$, so $\Sigma(E_a)=\{a\}$. Proposition~\ref{cl:essential-spectrum} distinguishes continuum many standard orbits, and Theorem~\ref{cl:main} distinguishes their metric types. Conversely, an operator on a separable Hilbert space is determined by countably many complex matrix entries, so $\mathcal B(H)$ has cardinality at most $2^{\aleph_0}$. This supplies the matching upper bound.
\end{proof}

The independent implementations are essential. If $V$ is a nonscalar unitary, the map that is the identity on $\mathcal F\cup\mathcal J$ and is $E\mapsto VEV^*$ on $\mathcal F^\vee$ is a surjective isometry for the extended $L^p$ metric. It cannot be a single standard transformation: fixing all finite-rank idempotents forces such a standard action to be the identity, whereas conjugation by $V$ acts nontrivially on $\mathcal F^\vee$. Thus a complete finite-distance class has rigid $L^p$ isometries, with independent choices across the extended metric partition.

\section{\texorpdfstring{Comparison of the $L^p$-metric regimes}{Comparison of the Lp-metric regimes}}\label{sec:comparison}
The classifications have the same spatial building blocks, but their domains impose different coherence conditions. In the table, $p=\infty$ refers only to the operator norm. The perturbation-class and extended-metric assertions concern $1\leq p<\infty$.

\begin{center}
\small
\renewcommand{\arraystretch}{1.2}
\begin{tabularx}{\textwidth}{@{}p{.24\textwidth}p{.26\textwidth}X@{}}
\toprule
Domain & Metric and hypotheses & Form of a surjective isometry\\
\midrule
$\Id_k(H)$, $1\leq k<\infty$ & $d_p$, $1\leq p\leq\infty$; $H$ infinite dimensional & One unitary or antiunitary implementation, with optional adjunction, on the whole stratum.\\
$\Id_{\infty,\infty}(H)$ & $d_\infty$; $H$ separable & One of the four standard forms throughout the sector.\\
$\mathcal C_p(E_0)$, infinite rank and corank & $d_p$, $1\leq p<\infty$; $H$ separable & One of the four standard forms on the entire $L^p$-perturbation class, including all its index components.\\
$\Id(H)$ & Extended $d_p$, $1\leq p<\infty$; $H$ separable & Independent implementations on finite-type sectors and a bijective matching of isometric infinite-rank finite-distance classes.\\
\bottomrule
\end{tabularx}
\end{center}

\begin{corollary}[Finite-corank version]\label{comparison:corank}
Let $1\leq k<\infty$ and equip
\[
\{E\in\Id(H):\rank(I-E)=k\}
\]
with the metric $d_p$, where $1\leq p\leq\infty$. Every surjective $L^p$-isometry is implemented by a unitary or antiunitary operator, possibly composed with adjunction.
\end{corollary}
\begin{proof}
The complement map $E\mapsto I-E$ is an isometry onto $\Id_k(H)$. Conjugating the fixed-rank classification by this map gives the claim.
\end{proof}

The tangent-space argument identifies projection points within a rank stratum without assuming a derivative of the original isometry. The relative-index argument couples the connected components inside one $L^p$-perturbation class and identifies exact orthogonal increments. Both arguments eventually produce finite-dimensional probes, but each retains the endpoint analysis required by its own metric geometry.

\vspace{0.2cm}
{\bf AI Statement:}  During the preparation of this manuscript, the author used GPT 6 Sol for searching necessary references, examining mathematical arguments, and language editing. The tangent geometry technique was also suggested by GPT 6 Sol to prove the projection preservation. The authors have checked and substantially revised all mathematical arguments generated by AI and take full responsibility for the final version.

\bigskip
\noindent 1.\ School of Mathematical Sciences, Chongqing Normal University, Chongqing, 401331, China

\smallskip
\noindent\textit{E-mail addresses:} \href{mailto:whqian86@163.com}{\texttt{whqian86@163.com}};\quad
\href{mailto:ererenyuze@163.com}{\texttt{ererenyuze@163.com}};\quad
\href{mailto:wuwm@amss.ac.cn}{\texttt{wuwm@amss.ac.cn}}
\end{document}